\documentclass[12pt,twoside]{amsart}
\usepackage[latin1]{inputenc}
\usepackage{amsmath, amsthm, amscd, amsfonts, amssymb, graphicx}
\usepackage[bookmarksnumbered, plainpages]{hyperref}

\newtheorem{thm}{Theorem}[section]

\newtheorem{lem}[thm]{Lemma}

\newtheorem{defn}[thm]{Definition}

\numberwithin{equation}{section}

\begin{document}

\title{\bf The $J$-twist $D_{J}$ of the Dirac operator and the Kastler-Kalau-Walze type theorem for six-dimensional manifolds with boundary}
\author{Siyao Liu \hskip 0.4 true cm  Yong Wang$^{*}$}

\thanks{{\scriptsize
\hskip -0.4 true cm \textit{2010 Mathematics Subject Classification:}
53C40; 53C42.
\newline \textit{$^{*}$Corresponding author}}}

\maketitle

\begin{abstract}
 \indent In \cite{LW2}, the authors proved a Kastler-Kalau-Walze type theorem for the $J$-twist $D_{J}$ of the Dirac operator on $3$-dimensional and $4$-dimensional almost product Riemannian spin manifold with boundary. In this paper, we develop the Kastler-Kalau-Walze type theorem for the $J$-twist $D_{J}$ of the Dirac operator on a $6$-dimensional almost product Riemannian spin manifold with boundary.
\end{abstract}

\vskip 0.2 true cm


\pagestyle{myheadings}
\markboth{\rightline {\scriptsize Liu}}
         {\leftline{\scriptsize The $J$-twist $D_{J}$ of the Dirac operator and the Kastler-Kalau-Walze type theorem}}

\bigskip
\bigskip


\section{ Introduction}
Based on the noncommutative residue found in \cite{Gu,Wo}, Connes coined the Kastler-Kalau-Walze theorem in $1995,$ which had been studied extensively by geometers \cite{Co1,Co2,Ka,KW,Ac,U}.
Wang generalized some results to the case of manifolds with boundary in \cite{Wa1,Wa2} and proved the Kastler-Kalau-Walze type theorem for the Dirac operator and the signature operator on lower-dimensional manifolds with boundary.
Most of the operators which have the leading symbol $\sqrt{-1}c(\xi),$ regarding the Kastler-Kalau-Walze theorem, have been studied extensively by, among others, the second author and by previous researchers \cite{Wa3,Wa4,Wa5,WW,WWW,Wa,LW1}.
However, Wu and Wang studied operators with the leading symbol $-\widehat{c}(V)c(\xi)$. In \cite{WW2}, they gave the proof of Kastler-Kalau-Walze type theorems  of the operators $\sqrt{-1}\widehat{c}(V)(d+\delta)$ and $-\sqrt{-1}\widehat{c}(V)(d+\delta)$ on $3, 4$-dimensional oriented compact manifolds with or without boundary.

On the other hand, Kim had given some preliminaries and lemmas about the Dirac operator $D$ and the $J$-twist in \cite{K}. In \cite{Chen1,Chen2}, the author got estimates on the higher eigenvalues of the Dirac operator on locally reducible Riemannian manifolds by the $J$-twist of the Dirac operator. By simple calculations, the leading symbol of the $J$-twist $D_{J}$ of the Dirac operator is not $\sqrt{-1}c(\xi)$.

To bring together two areas, Liu and Wang proved a Kastler-Kalau-Walze type theorem for the $J$-twist $D_{J}$ of the Dirac operator on $3$-dimensional and $4$-dimensional almost product Riemannian spin manifold with boundary in \cite{LW2}.
Unlike the case of the Dirac operator in \cite{Wa3}, the author found the boundary term is not zero for the $J$-twist $D_{J}$ of the Dirac operator.
In this paper, following the study \cite{LW2}, we develop the Kastler-Kalau-Walze type theorem for the $J$-twist $D_{J}$ of the Dirac operator on a $6$-dimensional almost product Riemannian spin manifold with boundary.

This paper is organized as follows.
In Section $2$, we give a brief exposition of the almost product Riemannian manifold and the $J$-twist $D_{J}$ of the Dirac operator.
Using the results in Sec.$2,$ we prove the Kastler-Kalau-Walze type theorem on $6$-dimensional almost product Riemannian spin manifold with boundary for the $J$-twist $D_{J}$ of the Dirac operator in the next section.


\vskip 1 true cm

\section{ The $J$-twist $D_{J}$ of the Dirac operator $D$ }

We give some definitions and basic notions which we will use in this paper.

Let $M$ be a $n$-dimensional ($n\geq 3$) oriented compact Riemannian manifold with a Riemannian metric $g^{M}$.
And let $\nabla^L$ be the Levi-Civita connection about $g^M$.
Suppose that $\partial_{i}$ is a natural local frame on $TM$ and $(g^{ij})_{1\leq i,j\leq n}$ is the inverse matrix associated to the metric matrix  $(g_{ij})_{1\leq i,j\leq n}$ on $M$.
We recall that the Dirac operator $D$ is locally given as following:
\begin{equation}
D=\sum_{i, j=1}^{n}g^{ij}c(\partial_{i})\nabla_{\partial_{j}}^{S}=\sum_{i=1}^{n}c(e_{i})\nabla_{e_{i}}^{S},
\end{equation}
where $c(e_{i})$ be the Clifford action which satisfies the relation
\begin{align}
&c(e_{i})c(e_{j})+c(e_{j})c(e_{i})=-2g^{M}(e_{i}, e_{j})=-2\delta_i^j,
\end{align}
\begin{align}
&\nabla_{\partial_{j}}^{S}=\partial_{i}+\sigma_{i}
\end{align}
and
\begin{align}
&\sigma_{i}=\frac{1}{4}\sum_{j, k=1}^{n}\langle \nabla_{\partial_{i}}^{L}e_{j}, e_{k}\rangle c(e_{j})c(e_{k}).
\end{align}

Let $J$ be a $(1, 1)$-tensor field on $(M, g^M)$ such that $J^2=\texttt{id},$
\begin{align}
&g^M(J(X), J(Y))=g^M(X, Y),
\end{align}
for all vector fields $X,Y\in \Gamma(TM).$ Here $\texttt{id}$ stands for the identity map. $(M, g^M, J)$ is an almost product Riemannian manifold. We can also define on almost product Riemannian spin manifold the following $J$-twist $D_{J}$ of the Dirac operator $D$ by
\begin{align}
&D_{J}:=\sum_{i=1}^{n}c(e_{i})\nabla^{S}_{J(e_{i})}=\sum_{i=1}^{n}c[J(e_{i})]\nabla^{S}_{e_{i}}.
\end{align}

Let $\xi=\sum_{k}\xi_{j}dx_{j},$  $\nabla^L_{\partial_{i}}\partial_{j}=\sum_{k}\Gamma_{ij}^{k}\partial_{k},$  we denote that
\begin{align}
 \Gamma^{k}=g^{ij}\Gamma_{ij}^{k};\ \sigma^{j}=g^{ij}\sigma_{i};\ \partial^{j}=g^{ij}\partial_{i},\nonumber
\end{align}
where $\Gamma_{ij}^{k}$ is the Christoffel coefficient of $\nabla^{L}$.

To shorten the notation, we let $s$ stand for the scalar curvature.
By (2.9) in \cite{LW2}, we have
\begin{align}
{D_{J}}^{2}&=-\frac{1}{8}\sum_{i,j,k,l=1}^{n}R(J(e_{i}), J(e_{j}), e_{k}, e_{l})c(e_{i})c(e_{j})c(e_{k})c(e_{l})-g^{ij}\partial_{i}\partial_{j}-2\sigma^{j}\partial_{j}+\Gamma^{k}\partial_{k}\\
&-g^{ij}[(\partial_{i}\sigma_{j})+\sigma_{i}\sigma_{j}-\Gamma^{k}_{ij}\sigma_{k}]+\frac{1}{4}s+\sum_{\alpha,\beta=1}^{n}c[J(e_{\alpha})]c[(\nabla^{L}_{e_{\alpha}}J)e_{\beta}]\sum_{k=1}^{n}\langle e_{\beta}, dx_{k}\rangle\nabla^{S}_{\partial_{k}},\nonumber
\end{align}
where $e_i^*=g^{M}(e_i,\cdot)$ and $\langle X, dx_{k}\rangle=g^{M}(X, (dx_{k})^{*}),$ for a vector field $X.$

Let us now turn to compute the specification of ${D_{J}}^{3}.$
\begin{align}
{D_{J}}^{3}&=\sum^{n}_{i=1}c[J(e_{i})]\langle e_{i},dx_{l}\rangle(-g^{ij}\partial_{l}\partial_{i}\partial_{j})+\sum^{n}_{i=1}c[J(e_{i})]\langle e_{i},dx_{l}\rangle\Big(-(\partial_{l}g^{ij})\partial_{i}\partial_{j}-g^{ij}(4\sigma_{i}\partial_{j}\\
&-2\Gamma^{k}_{ij}\partial_{k})\partial_{l}+2\sum^{n}_{\alpha,\beta,k=1}c[J(e_{\alpha})]c[(\nabla^{L}_{e_{\alpha}}J)e_{\beta}]\langle e_{\beta}, dx_{k}\rangle\partial_{l}\partial_{k}\Big)+\sum^{n}_{i=1}c[J(e_{i})]\sigma_{i}(-g^{ij}\partial_{i}\partial_{j})\nonumber\\
&+\sum^{n}_{i=1}c[J(e_{i})]\langle e_{i},dx_{l}\rangle\Big[-2(\partial_{l}g^{ij})\sigma_{i}\partial_{j}-2g^{ij}(\partial_{l}\sigma_{i})\partial_{j}+(\partial_{l}g^{ij})\Gamma^{k}_{ij}\partial_{k}+g^{ij}(\partial_{l}\Gamma^{k}_{ij})\partial_{k}\nonumber\\
&+\sum^{n}_{\alpha,\beta,k=1}\partial_{l}\Big(c[J(e_{\alpha})]c[(\nabla^{L}_{e_{\alpha}}J)e_{\beta}]\Big)\langle e_{\beta}, dx_{k}\rangle\partial_{k}+\sum^{n}_{\alpha,\beta,k=1}c[J(e_{\alpha})]c[(\nabla^{L}_{e_{\alpha}}J)e_{\beta}]\Big(\partial_{l}\langle e_{\beta}, dx_{k}\rangle\Big)\partial_{k}\Big]\nonumber\\
&+\sum^{n}_{i=1}c[J(e_{i})]\langle e_{i},dx_{l}\rangle\partial_{l}\Big(-\frac{1}{8}\sum_{i,j,k,l=1}^{n}R(J(e_{i}), J(e_{j}), e_{k}, e_{l})c(e_{i})c(e_{j})c(e_{k})c(e_{l})
-g^{ij}((\partial_{i}\sigma_{j})\nonumber\\
&+\sigma_{i}\sigma_{j}-\Gamma^{k}_{ij}\sigma_{k})+\frac{1}{4}s+\sum^{n}_{\alpha,\beta,k=1}c[J(e_{\alpha})]c[(\nabla^{L}_{e_{\alpha}}J)e_{\beta}]\langle e_{\beta}, dx_{k}\rangle\sigma_{k}\Big)+\sum^{n}_{i=1}c[J(e_{i})]\sigma_{i}\Big(\frac{1}{4}s\nonumber\\
&-\frac{1}{8}\sum_{i,j,k,l=1}^{n}R(J(e_{i}), J(e_{j}), e_{k}, e_{l})c(e_{i})c(e_{j})c(e_{k})c(e_{l})-2\sigma^{j}\partial_{j}+\Gamma^{k}\partial_{k}-g^{ij}((\partial_{i}\sigma_{j})+\sigma_{i}\sigma_{j}\nonumber\\
&-\Gamma^{k}_{ij}\sigma_{k})+\sum^{n}_{\alpha,\beta,k=1}c[J(e_{\alpha})]c[(\nabla^{L}_{e_{\alpha}}J)e_{\beta}]\langle e_{\beta}, dx_{k}\rangle(\partial_{k}+\sigma_{k})\Big).\nonumber
\end{align}

From \cite{LW2}, we also obtain
\begin{thm}\cite{LW2} If $M$ is a $n$-dimensional almost product Riemannian spin manifold without boundary, we have the following:
\begin{align}
{\rm Wres}({D_{J}}^{-n+2})
=\frac{(n-2)(4\pi)^{\frac{n}{2}}}{(\frac{n}{2}-1)!}\int_{M}2^{\frac{n}{2}}\Big(&\frac{1}{4}\sum_{i,j=1}^{n}R(J(e_{i}), J(e_{j}), e_{j}, e_{i})
-\frac{1}{2}\sum_{\nu,j=1}^{n}g^{M}(\nabla_{e_{j}}^{L}(J)e_{\nu}, (\nabla^{L}_{e_{\nu}}J)e_{j})\\
&-\frac{1}{2}\sum_{\nu,j=1}^{n}g^{M}(J(e_{\nu}), (\nabla^{L}_{e_{j}}(\nabla^{L}_{e_{\nu}}(J)))e_{j}-(\nabla^{L}_{\nabla^{L}_{e_{j}}e_{\nu}}(J))e_{j})\nonumber\\
&-\frac{1}{4}\sum_{\alpha,\nu,j=1}^{n}g^{M}(J(e_{\alpha}), (\nabla^{L}_{e_{\nu}}J)e_{j})g^{M}((\nabla^{L}_{e_{\alpha}}J)e_{j}, J(e_{\nu}))\nonumber\\
&-\frac{1}{4}\sum_{\alpha,\nu,j=1}^{n}g^{M}(J(e_{\alpha}), (\nabla^{L}_{e_{\alpha}}J)e_{j})g^{M}(J(e_{\nu}), (\nabla^{L}_{e_{\nu}}J)e_{j})\nonumber\\
&+\frac{1}{4}\sum_{\nu,j=1}^{n}g^{M}((\nabla^{L}_{e_{\nu}}J)e_{j}, (\nabla^{L}_{e_{\nu}}J)e_{j}))-\frac{1}{12}s\Big)d{\rm Vol_{M} }.\nonumber
\end{align}
\end{thm}

\section{ The Kastler-Kalau-Walze type theorem for $6$-dimensional manifolds with boundary }
\indent Firstly, we explain the basic notions of  Boutet de Monvel's calculus and the definition of the noncommutative residue for manifolds with boundary that will be used throughout the paper. For the details, see Ref.\cite{Wa3}.\\
\indent Let $U\subset M$ be a collar neighborhood of $\partial M$ which is diffeomorphic with $\partial M\times [0,1)$. By the definition of $h(x_n)\in C^{\infty}([0,1))$
and $h(x_n)>0$, there exists $\widehat{h}\in C^{\infty}((-\varepsilon,1))$ such that $\widehat{h}|_{[0,1)}=h$ and $\widehat{h}>0$ for some
sufficiently small $\varepsilon>0$. Then there exists a metric $g'$ on $\widetilde{M}=M\bigcup_{\partial M}\partial M\times
(-\varepsilon,0]$ which has the form on $U\bigcup_{\partial M}\partial M\times (-\varepsilon,0 ]$
\begin{equation}
g'=\frac{1}{\widehat{h}(x_{n})}g^{\partial M}+dx _{n}^{2} ,
\end{equation}
such that $g'|_{M}=g$. We fix a metric $g'$ on the $\widetilde{M}$ such that $g'|_{M}=g$.

We define the Fourier transformation $F'$  by
\begin{equation}
F':L^2({\bf R}_t)\rightarrow L^2({\bf R}_v);~F'(u)(v)=\int e^{-ivt}u(t)dt\\
\end{equation}
and let
\begin{equation}
r^{+}:C^\infty ({\bf R})\rightarrow C^\infty (\widetilde{{\bf R}^+});~ f\rightarrow f|\widetilde{{\bf R}^+};~
\widetilde{{\bf R}^+}=\{x\geq0;x\in {\bf R}\}.
\end{equation}
 where $\Phi({\bf R})$
denotes the Schwartz space and $\Phi(\widetilde{{\bf R}^+}) =r^+\Phi({\bf R})$, $\Phi(\widetilde{{\bf R}^-}) =r^-\Phi({\bf R})$. We define $H^+=F'(\Phi(\widetilde{{\bf R}^+}));~ H^-_0=F'(\Phi(\widetilde{{\bf R}^-}))$ which satisfies
$H^+\bot H^-_0$. We have the following
 property: $h\in H^+~(H^-_0)$ if and only if $h\in C^\infty({\bf R})$ which has an analytic extension to the lower (upper) complex
half-plane $\{{\rm Im}\xi<0\}~(\{{\rm Im}\xi>0\})$ such that for all nonnegative integer $l$,
 \begin{equation}
\frac{d^{l}h}{d\xi^l}(\xi)\sim\sum^{\infty}_{k=1}\frac{d^l}{d\xi^l}(\frac{c_k}{\xi^k}),
\end{equation}
as $|\xi|\rightarrow +\infty,{\rm Im}\xi\leq0~({\rm Im}\xi\geq0)$.\\
\indent Let $H'$ be the space of all polynomials and $H^-=H^-_0\bigoplus H';~H=H^+\bigoplus H^-.$ Denote by $\pi^+~(\pi^-)$ respectively the projection on $H^+~(H^-)$. For calculations, we take $H=\widetilde H=\{$rational functions having no poles on the real axis$\}$ ($\tilde{H}$ is a dense set in the topology of $H$). Then on $\tilde{H}$,
 \begin{equation}
\pi^+h(\xi_0)=\frac{1}{2\pi i}\lim_{u\rightarrow 0^{-}}\int_{\Gamma^+}\frac{h(\xi)}{\xi_0+iu-\xi}d\xi,
\end{equation}
where $\Gamma^+$ is a Jordan close curve
included ${\rm Im}(\xi)>0$ surrounding all the singularities of $h$ in the upper half-plane and
$\xi_0\in {\bf R}$. Similarly, define $\pi'$ on $\tilde{H}$,
\begin{equation}
\pi'h=\frac{1}{2\pi}\int_{\Gamma^+}h(\xi)d\xi.
\end{equation}
So, $\pi'(H^-)=0$. For $h\in H\bigcap L^1({\bf R})$, $\pi'h=\frac{1}{2\pi}\int_{{\bf R}}h(v)dv$ and for $h\in H^+\bigcap L^1({\bf R})$, $\pi'h=0$.

Let $M$ be a $n$-dimensional compact oriented manifold with boundary $\partial M$.
Denote by $\mathcal{B}$ Boutet de Monvel's algebra, we recall the main theorem in \cite{Wa3,FGLS}.
\begin{thm}\label{th:32}{\rm\cite{FGLS}}{\bf(Fedosov-Golse-Leichtnam-Schrohe)}
 Let $X$ and $\partial X$ be connected, ${\rm dim}X=n\geq3$,
 $A=\left(\begin{array}{lcr}\pi^+P+G &   K \\
T &  S    \end{array}\right)$ $\in \mathcal{B}$ , and denote by $p$, $b$ and $s$ the local symbols of $P,G$ and $S$ respectively.
 Define:
 \begin{align}
{\rm{\widetilde{Wres}}}(A)&=\int_X\int_{\bf S}{\rm{tr}}_E\left[p_{-n}(x,\xi)\right]\sigma(\xi)dx \\
&+2\pi\int_ {\partial X}\int_{\bf S'}\left\{{\rm tr}_E\left[({\rm{tr}}b_{-n})(x',\xi')\right]+{\rm{tr}}
_F\left[s_{1-n}(x',\xi')\right]\right\}\sigma(\xi')dx'.\nonumber
\end{align}
Then~~ a) ${\rm \widetilde{Wres}}([A,B])=0 $, for any
$A,B\in\mathcal{B}$;~~ b) It is a unique continuous trace on
$\mathcal{B}/\mathcal{B}^{-\infty}$.
\end{thm}

\begin{defn}{\rm\cite{Wa3} }
Lower dimensional volumes of spin manifolds with boundary are defined by
 \begin{equation}
{\rm Vol}^{(p_1,p_2)}_nM:= \widetilde{{\rm Wres}}[\pi^+D^{-p_1}\circ\pi^+D^{-p_2}].
\end{equation}
\end{defn}
We can get the spin structure on $\widetilde{M}$ by extending the spin structure on $M.$ Let $D$ be the Dirac operator associated to $g'$ on the spinors bundle $S(T\widetilde{M}).$ By \cite{Wa3}, we get
\begin{align}
\widetilde{{\rm Wres}}[\pi^+D^{-p_1}\circ\pi^+D^{-p_2}]=\int_M\int_{|\xi|=1}{\rm
trace}_{S(TM)}[\sigma_{-n}(D^{-p_1-p_2})]\sigma(\xi)dx+\int_{\partial M}\Phi
\end{align}
and
\begin{align}
\Phi&=\int_{|\xi'|=1}\int^{+\infty}_{-\infty}\sum^{\infty}_{j, k=0}\sum\frac{(-i)^{|\alpha|+j+k+1}}{\alpha!(j+k+1)!}
\times {\rm trace}_{S(TM)}[\partial^j_{x_n}\partial^\alpha_{\xi'}\partial^k_{\xi_n}\sigma^+_{r}(D^{-p_1})(x',0,\xi',\xi_n)
\\
&\times\partial^\alpha_{x'}\partial^{j+1}_{\xi_n}\partial^k_{x_n}\sigma_{l}(D^{-p_2})(x',0,\xi',\xi_n)]d\xi_n\sigma(\xi')dx',\nonumber
\end{align}
 where the sum is taken over $r+l-k-|\alpha|-j-1=-n,~~r\leq -p_1,l\leq -p_2$.

 Since $[\sigma_{-n}(D^{-p_1-p_2})]|_M$ has the same expression as $\sigma_{-n}(D^{-p_1-p_2})$ in the case of manifolds without
boundary, so locally we can compute the first term by \cite{Ka}, \cite{KW}, \cite{Wa3}, \cite{Po}.

For any fixed point $x_0\in\partial M$, we choose the normal coordinates
$U$ of $x_0$ in $\partial M$ (not in $M$) and compute $\Phi(x_0)$ in the coordinates $\widetilde{U}=U\times [0,1)\subset M$ and the
metric $\frac{1}{h(x_n)}g^{\partial M}+dx_n^2.$ The dual metric of $g^M$ on $\widetilde{U}$ is ${h(x_n)}g^{\partial M}+dx_n^2.$  Write
$g^M_{ij}=g^M(\frac{\partial}{\partial x_i},\frac{\partial}{\partial x_j});~ g_M^{ij}=g^M(dx_i,dx_j)$, then

\begin{equation}
[g^M_{i,j}]= \left[\begin{array}{lcr}
  \frac{1}{h(x_n)}[g_{i,j}^{\partial M}]  & 0  \\
   0  &  1
\end{array}\right];~~~
[g_M^{i,j}]= \left[\begin{array}{lcr}
  h(x_n)[g^{i,j}_{\partial M}]  & 0  \\
   0  &  1
\end{array}\right]
\end{equation}
and
\begin{equation}
\partial_{x_s}g_{ij}^{\partial M}(x_0)=0, 1\leq i,j\leq n-1; ~~~g_{ij}^M(x_0)=\delta_{ij}.
\end{equation}
\indent $\{e_1, \cdots, e_n\}$ be an orthonormal frame field in $U$ about $g^{\partial M}$ which is parallel along geodesics and $e_i(x_0)=\frac{\partial}{\partial{x_i}}(x_0).$ We review the following three lemmas.
\begin{lem}{\rm \cite{Wa3}}\label{le:32}
With the metric $g^{M}$ on $M$ near the boundary
\begin{eqnarray}
\partial_{x_j}(|\xi|_{g^M}^2)(x_0)&=&\left\{
       \begin{array}{c}
        0,  ~~~~~~~~~~ ~~~~~~~~~~ ~~~~~~~~~~~~~{\rm if }~j<n, \\[2pt]
       h'(0)|\xi'|^{2}_{g^{\partial M}},~~~~~~~~~~~~~~~~~~~~{\rm if }~j=n;
       \end{array}
    \right. \\
\partial_{x_j}[c(\xi)](x_0)&=&\left\{
       \begin{array}{c}
      0,  ~~~~~~~~~~ ~~~~~~~~~~ ~~~~~~~~~~~~~{\rm if }~j<n,\\[2pt]
\partial x_{n}(c(\xi'))(x_{0}), ~~~~~~~~~~~~~~~~~{\rm if }~j=n,
       \end{array}
    \right.
\end{eqnarray}
where $\xi=\xi'+\xi_{n}dx_{n}$.
\end{lem}
\begin{lem}{\rm \cite{Wa3}}\label{le:32}With the metric $g^{M}$ on $M$ near the boundary
\begin{align}
\omega_{s,t}(e_i)(x_0)&=\left\{
       \begin{array}{c}
        \omega_{n,i}(e_i)(x_0)=\frac{1}{2}h'(0),  ~~~~~~~~~~ ~~~~~~~~~~~{\rm if }~s=n,t=i,i<n; \\[2pt]
       \omega_{i,n}(e_i)(x_0)=-\frac{1}{2}h'(0),~~~~~~~~~~~~~~~~~~~{\rm if }~s=i,t=n,i<n;\\[2pt]
    \omega_{s,t}(e_i)(x_0)=0,~~~~~~~~~~~~~~~~~~~~~~~~~~~other~cases,~~~~~~~~~\\[2pt]
       \end{array}
    \right.
\end{align}
where $(\omega_{s,t})$ denotes the connection matrix of Levi-Civita connection $\nabla^L$.
\end{lem}
\begin{lem}{\rm \cite{Wa3}}
\begin{align}
\Gamma_{st}^k(x_0)&=\left\{
       \begin{array}{c}
        \Gamma^n_{ii}(x_0)=\frac{1}{2}h'(0),~~~~~~~~~~ ~~~~~~~~~~~{\rm if }~s=t=i,k=n,i<n; \\[2pt]
        \Gamma^i_{ni}(x_0)=-\frac{1}{2}h'(0),~~~~~~~~~~~~~~~~~~~{\rm if }~s=n,t=i,k=i,i<n;\\[2pt]
        \Gamma^i_{in}(x_0)=-\frac{1}{2}h'(0),~~~~~~~~~~~~~~~~~~~{\rm if }~s=i,t=n,k=i,i<n,\\[2pt]
        \Gamma_{st}^i(x_0)=0,~~~~~~~~~~~~~~~~~~~~~~~~~~~other~cases.~~~~~~~~~
       \end{array}
    \right.
\end{align}
\end{lem}
\indent Applying (3.9) and (3.10) yields
\begin{equation}
\widetilde{{\rm Wres}}[\pi^+{{D}_{J}}^{-1}\circ\pi^+{{D}_{J}}^{-3}]=\int_M\int_{|\xi|=1}{\rm
trace}_{S(TM)}[\sigma_{-6}({{D}_{J}}^{-4})]\sigma(\xi)dx+\int_{\partial M}\Phi,
\end{equation}
where
\begin{align}
\Phi &=\int_{|\xi'|=1}\int^{+\infty}_{-\infty}\sum^{\infty}_{j, k=0}\sum\frac{(-i)^{|\alpha|+j+k+1}}{\alpha!(j+k+1)!}
\times {\rm trace}_{S(TM)}[\partial^j_{x_n}\partial^\alpha_{\xi'}\partial^k_{\xi_n}\sigma^+_{r}({{D}_{J}}^{-1})\\
&(x',0,\xi',\xi_n)\times\partial^\alpha_{x'}\partial^{j+1}_{\xi_n}\partial^k_{x_n}\sigma_{l}({{D}_{J}}^{-3})(x',0,\xi',\xi_n)]d\xi_n\sigma(\xi')dx',\nonumber
\end{align}
the sum is taken over $r+l-k-j-|\alpha|-1=-6, r\leq -1, l\leq-3$.\\

\indent
Computations show that
\begin{align}
&\int_M\int_{|\xi|=1}{\rm trace}_{S(TM)}[\sigma_{-6}({{D}_{J}}^{-4})]\sigma(\xi)dx=256\pi^{3}\\
&\int_{M}\Big(\sum_{i,j=1}^{n}R(J(e_{i}), J(e_{j}), e_{j}, e_{i})
-2\sum_{\nu,j=1}^{n}g^{M}(\nabla_{e_{j}}^{L}(J)e_{\nu}, (\nabla^{L}_{e_{\nu}}J)e_{j})\nonumber\\
&-2\sum_{\nu,j=1}^{n}g^{M}(J(e_{\nu}), (\nabla^{L}_{e_{j}}(\nabla^{L}_{e_{\nu}}(J)))e_{j}-(\nabla^{L}_{\nabla^{L}_{e_{j}}e_{\nu}}(J))e_{j})\nonumber\\
&-\sum_{\alpha,\nu,j=1}^{n}g^{M}(J(e_{\alpha}), (\nabla^{L}_{e_{\nu}}J)e_{j})g^{M}((\nabla^{L}_{e_{\alpha}}J)e_{j}, J(e_{\nu}))\nonumber\\
&-\sum_{\alpha,\nu,j=1}^{n}g^{M}(J(e_{\alpha}), (\nabla^{L}_{e_{\alpha}}J)e_{j})g^{M}(J(e_{\nu}), (\nabla^{L}_{e_{\nu}}J)e_{j})\nonumber\\
&+\sum_{\nu,j=1}^{n}g^{M}((\nabla^{L}_{e_{\nu}}J)e_{j}, (\nabla^{L}_{e_{\nu}}J)e_{j}))-\frac{1}{3}s\Big)d{\rm Vol_{M} },\nonumber
\end{align}
where $\Omega_{n}=\frac{2\pi^\frac{n}{2}}{\Gamma(\frac{n}{2})}.$

\indent Now, we compute $\int_{\partial M} \Phi.$
As shown in \cite{LW2}, we see that
\begin{lem}\cite{LW2} The following identities hold:
\begin{align}
\sigma_1({D}_{J})&=ic[J(\xi)];\\
\sigma_0({D}_{J})&=
-\frac{1}{4}\sum_{i,j,k=1}^{n}\omega_{j,k}(e_i)c[J(e_i)]c(e_j)c(e_k).
\end{align}
\end{lem}
\begin{lem}\cite{LW2} The following identities hold:
\begin{align}
\sigma_{-1}({{D}_{J}}^{-1})&=\frac{ic[J(\xi)]}{|\xi|^2};\\
\sigma_{-2}({{D}_{J}}^{-1})&=\frac{c[J(\xi)]\sigma_{0}({D}_{J})c[J(\xi)]}{|\xi|^4}+\frac{c[J(\xi)]}{|\xi|^6}\sum_ {j=1}^{n} c[J(dx_j)]
\Big[\partial_{x_j}(c[J(\xi)])|\xi|^2-c[J(\xi)]\partial_{x_j}(|\xi|^2)\Big].
\end{align}
\end{lem}
According to (2.9), we have
\begin{lem} The following identities hold:
\begin{align}
\sigma_3({{D}_{J}}^{3})&=ic[J(\xi)]|\xi|^{2};\\
\sigma_2({{D}_{J}}^{3})&=\sum^{n}_{i,j,l=1}c[J(dx_{l})]\partial_{l}(g^{ij})\xi_{i}\xi_{j}+c[J(\xi)](4\sigma^k-2\Gamma^k)\xi_{k}-2\sum^{n}_{\alpha=1}c[J(\xi)]c[J(e_{\alpha})]c[(\nabla^{L}_{e_{\alpha}}J)(\xi^{*})]\\
&-\frac{1}{4}|\xi|^2\sum^{n}_{s,t,l=1}\omega_{s,t}(e_{l})c[J(e_{l})]c(e_{s})c(e_{t}),\nonumber
\end{align}
where $\xi^{*}=\sum^{n}_{\beta=1}\langle e_{\beta}, \xi\rangle e_{\beta}.$
\end{lem}
\indent Write
 \begin{eqnarray}
D_x^{\alpha}&=(-i)^{|\alpha|}\partial_x^{\alpha};
~\sigma({{D}_{J}}^{3})=p_{3}+p_{2}+p_{1}+p_{0};
~\sigma({{D}_{J}}^{-3})=\sum^{\infty}_{j=3}q_{-j}.
\end{eqnarray}

\indent By the composition formula of pseudodifferential operators, we have
\begin{align}
1=\sigma({{D}_{J}}^3\circ {{D}_{J}}^{-3})&=
\sum_{\alpha}\frac{1}{\alpha!}\partial^{\alpha}_{\xi}
[\sigma({{D}_{J}}^3)]{{D}}^{\alpha}_{x}
[\sigma({{D}_{J}}^{-3})] \\
&=(p_3+p_2+p_1+p_0)(q_{-3}+q_{-4}+q_{-5}+\cdots) \nonumber\\
&+\sum_j(\partial_{\xi_{j}}p_3+\partial_{\xi_{j}}p_2+\partial_{\xi_j}p_1+\partial_{\xi_{j}}p_0)
(D_{x_j}q_{-3}+D_{x_j}q_{-4}+D_{x_j}q_{-5}+\cdots) \nonumber\\
&=p_3q_{-3}+(p_3q_{-4}+p_2q_{-3}+\sum_j\partial_{\xi_j}p_3D_{x_j}q_{-3})+\cdots,\nonumber
\end{align}
so
\begin{equation}
q_{-3}=p_3^{-1};~q_{-4}=-p_3^{-1}[p_2p_3^{-1}+\sum_j\partial_{\xi_j}p_3D_{x_j}(p_3^{-1})].
\end{equation}
\begin{lem} The following identities hold:
\begin{align}
\sigma_{-3}({{D}_{J}}^{-3})&=\frac{ic[J(\xi)]}{|\xi|^{4}};\\
\sigma_{-4}({{D}_{J}}^{-3})&=\frac{c[J(\xi)]\sigma_{2}({{D}_{J}}^{3})c[J(\xi)]}{|\xi|^{8}}+\frac{c[J(\xi)]}{|\xi|^{10}}\sum_ {j=1}^{n} \Big(c[J(dx_j)]|\xi|^{2}+2\xi_{j}c[J(\xi)]\Big)\Big[\partial_{x_j}(c[J(\xi)])|\xi|^2\\
&-2c[J(\xi)]\partial_{x_j}(|\xi|^2)\Big].\nonumber
\end{align}
\end{lem}
\indent When $n=6$, then ${\rm tr}[{\rm \texttt{id}}]=8,$ since the sum is taken over $
r+l-k-j-|\alpha|-1=-6,~~r\leq -1,l\leq-3,$ then we have the following five cases:\\

\noindent  {\bf case (a)~(I)}~$r=-1, l=-3, j=k=0, |\alpha|=1$.\\

\noindent By applying the formula shown in (3.18), we can calculate
\begin{equation}
\Phi_1=-\int_{|\xi'|=1}\int^{+\infty}_{-\infty}\sum_{|\alpha|=1}{\rm trace}
[\partial^{\alpha}_{\xi'}\pi^{+}_{\xi_{n}}\sigma_{-1}({{D}_{J}}^{-1})\times\partial^{\alpha}_{x'}\partial_{\xi_{n}}\sigma_{-3}({{D}_{J}}^{-3})](x_0)d\xi_n\sigma(\xi')dx'.
\end{equation}
For $i<n,$ we get
\begin{align}
\partial_{x_i}\left(\frac{ic[J(\xi)]}{|\xi|^{4}}\right)(x_0)
=\frac{i\partial_{x_i}(c[J(\xi)])(x_0)}{|\xi|^{4}}-\frac{ic[J(\xi)]\partial_{x_i}(|\xi|^{4})(x_0)}{|\xi|^{8}}
=\frac{i\partial_{x_i}(c[J(\xi)])(x_0)}{|\xi|^{4}},
\end{align}
where $J(dx_{p})=\sum^{n}_{h=1}a^{p}_{h}dx_{h}.$\\
Of course,
\begin{align}
\partial_{x_i}\left(\frac{ic[J(\xi)]}{|\xi|^{4}}\right)(x_0)
&=\frac{i\partial_{x_i}(c[J(\sum^{n}_{p=1}\xi_{p}dx_{p})])(x_0)}{|\xi|^{4}}=\frac{i\sum^{n}_{p=1}\xi_{p}\partial_{x_i}(c(\sum^{n}_{h=1}a^{p}_{h}dx_{h}))(x_0)}{|\xi|^{4}}\\
&=\frac{i\sum^{n}_{p,h=1}\xi_{p}\partial_{x_i}(a^{p}_{h})c(dx_{h})(x_0)}{|\xi|^{4}}+\frac{i\sum^{n}_{p,h=1}\xi_{p}a^{p}_{h}\partial_{x_i}(c(dx_{h}))(x_0)}{|\xi|^{4}}\nonumber\\
&=\frac{i\sum^{n}_{p,h=1}\xi_{p}\partial_{x_i}(a^{p}_{h})c(dx_{h})(x_0)}{|\xi|^{4}}.\nonumber
\end{align}
When $|\xi'|=1$, we see that
\begin{align}
\partial_{x_i}\left(\frac{ic[J(\xi)]}{|\xi|^{4}}\right)(x_0)|_{|\xi'|=1}
&=\frac{i\sum^{n}_{h=1}\sum^{n-1}_{p=1}\xi_{p}\partial_{x_i}(a^{p}_{h})c(dx_{h})(x_0)}{(1+\xi_{n}^{2})^{2}}+\frac{i\sum^{n}_{h=1}\xi_{n}\partial_{x_i}(a^{n}_{h})c(dx_{h})(x_0)}{(1+\xi_{n}^{2})^{2}}.
\end{align}
An easy computation shows that
\begin{align}
\partial_{\xi_{n}}\partial_{x_i}\left(\frac{ic[J(\xi)]}{|\xi|^{4}}\right)(x_0)|_{|\xi'|=1}
&=i\sum^{n}_{h=1}\sum^{n-1}_{p=1}\xi_{p}\partial_{x_i}(a^{p}_{h})c(dx_{h})\partial_{\xi_{n}}\left(\frac{1}{(1+\xi_{n}^{2})^{2}}\right)\\
&+i\sum^{n}_{h=1}\partial_{x_i}(a^{n}_{h})c(dx_{h})\partial_{\xi_{n}}\left(\frac{\xi_{n}}{(1+\xi_{n}^{2})^{2}}\right)\nonumber\\
&=\frac{-4i\xi_{n}}{(1+\xi_{n}^{2})^{3}}\sum^{n}_{h=1}\sum^{n-1}_{p=1}\xi_{p}\partial_{x_i}(a^{p}_{h})c(dx_{h})\nonumber\\
&+\frac{i(1-3\xi_{n}^2)}{(1+\xi_{n}^{2})^{3}}\sum^{n}_{h=1}\partial_{x_i}(a^{n}_{h})c(dx_{h}).\nonumber
\end{align}
(3.39) in \cite{LW2} makes it obvious that
\begin{align}
\pi^+_{\xi_n}\partial_{\xi_i}\left(\frac{ic[J(\xi)]}{|\xi|^2}\right)(x_0)|_{|\xi'|=1}
&=\frac{1}{2(\xi_n-i)}c[J(dx_i)]+\frac{2i-\xi_n}{2(\xi_n-i)^2}\sum^{n-1}_{q=1}\xi_{i}\xi_{q}c[J(dx_q)]\\
&-\frac{1}{2(\xi_n-i)^2}\xi_{i}c[J(dx_n)].\nonumber
\end{align}
Substituting $c[J(dx_{p})]$ into $\sum^{n}_{h=1}a^{p}_{h}c(dx_{h}),$ we can rewrite (3.36) as
\begin{align}
\pi^+_{\xi_n}\partial_{\xi_i}\left(\frac{ic[J(\xi)]}{|\xi|^2}\right)(x_0)|_{|\xi'|=1}
&=\frac{1}{2(\xi_n-i)}\sum^{n}_{\alpha=1}a^{i}_{\alpha}c(dx_{\alpha})+\frac{2i-\xi_n}{2(\xi_n-i)^2}\sum^{n}_{\beta=1}\sum^{n-1}_{q=1}\xi_{i}\xi_{q}a^{q}_{\beta}c(dx_{\beta})\\
&-\frac{1}{2(\xi_n-i)^2}\sum^{n}_{\gamma=1}\xi_{i}a^{n}_{\gamma}c(dx_{\gamma}).\nonumber
\end{align}
Hence, we have
\begin{align}
&\sum_{|\alpha|=1}{\rm trace}[\partial^\alpha_{\xi'}\pi^+_{\xi_n}\sigma_{-1}({{D}_{J}}^{-1})\times\partial^\alpha_{x'}\partial_{\xi_n}\sigma_{-3}({{D}_{J}}^{-3})](x_0)|_{|\xi'|=1}\\
=&-\frac{2 i \xi _n}{\left(\xi _n-i\right)^4 \left(\xi _n+i\right)^3}\sum_{\alpha,h=1}^{n}\sum_{i,p=1}^{n-1}{\rm tr}[\xi_{p}a_{\alpha}^{i}\partial_{x_i}(a_{h}^{p})c(dx_{\alpha})c(dx_{h})]\nonumber\\
&+\frac{i \left(1-3 \xi _n^2\right)}{2 \left(\xi _n-i\right)^4 \left(\xi _n+i\right)^3}\sum_{\alpha,h=1}^{n}\sum_{i=1}^{n-1}{\rm tr}[a_{\alpha}^{i}\partial_{x_i}(a_{h}^{n})c(dx_{\alpha})c(dx_{h})]\nonumber\\
&+\frac{2 \left(2+i \xi _n\right) \xi _n}{\left(\xi _n-i\right)^5 \left(\xi _n+i\right)^3}\sum_{\beta,h=1}^{n}\sum_{i,q,p=1}^{n-1}{\rm tr}[\xi_{i}\xi_{q}\xi_{p}a_{\beta}^{q}\partial_{x_i}(a_{h}^{p})c(dx_{\beta})c(dx_{h})]\nonumber\\
&+\frac{\left(2+i \xi _n\right) \left(3 \xi _n^2-1\right)}{2 \left(\xi _n-i\right)^5 \left(\xi _n+i\right)^3}\sum_{\beta,h=1}^{n}\sum_{i,q=1}^{n-1}{\rm tr}[\xi_{i}\xi_{q}a_{\beta}^{q}\partial_{x_i}(a_{h}^{n})c(dx_{\beta})c(dx_{h})]\nonumber\\
&+\frac{2 i \xi _n}{\left(\xi _n-i\right)^5 \left(\xi _n+i\right)^3}\sum_{\gamma,h=1}^{n}\sum_{i,p=1}^{n-1}{\rm tr}[\xi_{i}\xi_{p}a_{\gamma}^{n}\partial_{x_i}(a_{h}^{p})c(dx_{\gamma})c(dx_{h})]\nonumber\\
&+\frac{i \left(3 \xi _n^2-1\right)}{2 \left(\xi _n-i\right)^5 \left(\xi _n+i\right)^3}\sum_{\gamma,h=1}^{n}\sum_{i=1}^{n-1}{\rm tr}[\xi_{i}a_{\gamma}^{n}\partial_{x_i}(a_{h}^{n})c(dx_{\gamma})c(dx_{h})].\nonumber
\end{align}
Because $c(e_i)c(e_j)+c(e_j)c(e_i)=-2\delta_i^j$ then by the relation of the Clifford action and ${\rm tr}{AB}={\rm tr}{BA}$, we have the following equalities:
\begin{align}
\sum_{\alpha,h=1}^{n}\sum_{i,p=1}^{n-1}{\rm tr}[\xi_{p}a_{\alpha}^{i}\partial_{x_i}(a_{h}^{p})c(dx_{\alpha})c(dx_{h})]=-\sum_{h=1}^{n}\sum_{i,p=1}^{n-1}\xi_{p}a_{h}^{i}\partial_{x_i}(a_{h}^{p}){\rm tr}[\texttt{id}];
\end{align}
\begin{align}
\sum_{\alpha,h=1}^{n}\sum_{i=1}^{n-1}{\rm tr}[a_{\alpha}^{i}\partial_{x_i}(a_{h}^{n})c(dx_{\alpha})c(dx_{h})]=-\sum_{h=1}^{n}\sum_{i=1}^{n-1}a_{h}^{i}\partial_{x_i}(a_{h}^{n}){\rm tr}[\texttt{id}];
\end{align}
\begin{align}
\sum_{\beta,h=1}^{n}\sum_{i,q,p=1}^{n-1}{\rm tr}[\xi_{i}\xi_{q}\xi_{p}a_{\beta}^{q}\partial_{x_i}(a_{h}^{p})c(dx_{\beta})c(dx_{h})]=-\sum_{h=1}^{n}\sum_{i,q,p=1}^{n-1}\xi_{i}\xi_{q}\xi_{p}a_{h}^{q}\partial_{x_i}(a_{h}^{p})
{\rm tr}[\texttt{id}];
\end{align}
\begin{align}
\sum_{\beta,h=1}^{n}\sum_{i,q=1}^{n-1}{\rm tr}[\xi_{i}\xi_{q}a_{\beta}^{q}\partial_{x_i}(a_{h}^{n})c(dx_{\beta})c(dx_{h})]=-\sum_{h=1}^{n}\sum_{i,q=1}^{n-1}\xi_{i}\xi_{q}a_{h}^{q}\partial_{x_i}(a_{h}^{n}){\rm tr}[\texttt{id}];
\end{align}
\begin{align}
\sum_{\gamma,h=1}^{n}\sum_{i,p=1}^{n-1}{\rm tr}[\xi_{i}\xi_{p}a_{\gamma}^{n}\partial_{x_i}(a_{h}^{p})c(dx_{\gamma})c(dx_{h})]=-\sum_{h=1}^{n}\sum_{i,p=1}^{n-1}\xi_{i}\xi_{p}a_{h}^{n}\partial_{x_i}(a_{h}^{p}){\rm tr}[\texttt{id}];
\end{align}
\begin{align}
\sum_{\gamma,h=1}^{n}\sum_{i=1}^{n-1}{\rm tr}[\xi_{i}a_{\gamma}^{n}\partial_{x_i}(a_{h}^{n})c(dx_{\gamma})c(dx_{h})]=-\sum_{h=1}^{n}\sum_{i=1}^{n-1}\xi_{i}a_{h}^{n}\partial_{x_i}(a_{h}^{n}){\rm tr}[\texttt{id}].
\end{align}
We note that $\int_{|\xi'|=1}{\{\xi_{i_1}\cdot\cdot\cdot\xi_{i_{2d+1}}}\}\sigma(\xi')=0,$ this gives
\begin{align}
\Phi_1&=-\int_{|\xi'|=1}\int^{+\infty}_{-\infty}\sum_{|\alpha|=1}{\rm trace}
[\partial^{\alpha}_{\xi'}\pi^{+}_{\xi_{n}}\sigma_{-1}({{D}_{J}}^{-1})\times\partial^{\alpha}_{x'}\partial_{\xi_{n}}\sigma_{-3}({{D}_{J}}^{-3})](x_0)d\xi_n\sigma(\xi')dx'\\
&=\int_{|\xi'|=1}\int^{+\infty}_{-\infty}\frac{i \left(1-3 \xi _n^2\right)}{2 \left(\xi _n-i\right)^4 \left(\xi _n+i\right)^3}\sum_{h=1}^{n}\sum_{i=1}^{n-1}a_{h}^{i}\partial_{x_i}(a_{h}^{n}){\rm tr}[\texttt{id}]d\xi_n\sigma(\xi')dx'\nonumber\\
&+\int_{|\xi'|=1}\int^{+\infty}_{-\infty}\frac{\left(2+i \xi _n\right) \left(3 \xi _n^2-1\right)}{2 \left(\xi _n-i\right)^5 \left(\xi _n+i\right)^3}\sum_{h=1}^{n}\sum_{i,q=1}^{n-1}\xi_{i}\xi_{q}a_{h}^{q}\partial_{x_i}(a_{h}^{n}){\rm tr}[\texttt{id}]d\xi_n\sigma(\xi')dx'\nonumber\\
&+\int_{|\xi'|=1}\int^{+\infty}_{-\infty}\frac{2 i \xi _n}{\left(\xi _n-i\right)^5 \left(\xi _n+i\right)^3}\sum_{h=1}^{n}\sum_{i,p=1}^{n-1}\xi_{i}\xi_{p}a_{h}^{n}\partial_{x_i}(a_{h}^{p}){\rm tr}[\texttt{id}]d\xi_n\sigma(\xi')dx'\nonumber.\nonumber
\end{align}
From \cite{Ka}, we obtain $\int_{|\xi'|=1}\xi_i\xi_j=\frac{8\pi^{2}}{15}\delta_i^j,$ then
\begin{align}
\Phi_1
&=\sum_{h=1}^{n}\sum_{i=1}^{n-1}a_{h}^{i}\partial_{x_i}(a_{h}^{n}){\rm tr}[\texttt{id}]\Omega_4\frac{2\pi i}{3!}\Big[\frac{i \left(1-3 \xi _n^2\right)}{2 \left(\xi _n+i\right)^3}\Big]^{(3)}\bigg|_{\xi_n=i}dx'\\
&+\sum_{h=1}^{n}\sum_{i=1}^{n-1}a_{h}^{i}\partial_{x_i}(a_{h}^{n}){\rm tr}[\texttt{id}]\Omega_4\frac{8\pi^{2}}{15}\frac{2\pi i}{4!}\Big[\frac{\left(2+i \xi _n\right) \left(3 \xi _n^2-1\right)}{2 \left(\xi _n+i\right)^3}\Big]^{(4)}\bigg|_{\xi_n=i}dx'\nonumber\\
&+\sum_{h=1}^{n}\sum_{i=1}^{n-1}a_{h}^{n}\partial_{x_i}(a_{h}^{i}){\rm tr}[\texttt{id}]\Omega_4\frac{8\pi^{2}}{15}\frac{2\pi i}{4!}\Big[\frac{2 i \xi _n}{\left(\xi _n+i\right)^3}\Big]^{(4)}\bigg|_{\xi_n=i}dx'\nonumber\\
&=\sum_{h=1}^{n}\sum_{i=1}^{n-1}a_{h}^{i}\partial_{x_i}(a_{h}^{n}){\rm tr}[\texttt{id}]\Omega_4(-\frac{\pi}{16}+\frac{\pi^{3}}{12})dx'+\sum_{h=1}^{n}\sum_{i=1}^{n-1}a_{h}^{n}\partial_{x_i}(a_{h}^{i}){\rm tr}[\texttt{id}]\Omega_4(-\frac{\pi^{3}}{12})dx'.\nonumber
\end{align}

\noindent  {\bf case (a)~(II)}~$r=-1, l=-3, |\alpha|=k=0, j=1$.\\

\noindent It is easy to check that
\begin{equation}
\Phi_2=-\frac{1}{2}\int_{|\xi'|=1}\int^{+\infty}_{-\infty} {\rm trace}[\partial_{x_{n}}\pi^{+}_{\xi_{n}}\sigma_{-1}({{D}_{J}}^{-1})\times\partial^{2}_{\xi_{n}}\sigma_{-3}({{D}_{J}}^{-3})](x_0)d\xi_n\sigma(\xi')dx'.
\end{equation}
We can assert that
\begin{align}
&\pi^+_{\xi_n}\partial_{x_n}\left(\frac{ic[J(\xi)]}{|\xi|^2}\right)(x_0)|_{|\xi'|=1}\\
&=\frac{1}{2(\xi_n-i)}\sum^{n}_{h=1}\sum^{n-1}_{p=1}\xi_{p}\partial_{x_n}(a^{p}_{h})c(dx_{h})+\frac{i}{2(\xi_n-i)}\sum^{n}_{h=1}\partial_{x_n}(a^{n}_{h})c(dx_{h})\nonumber\\
&+\frac{1}{2(\xi_n-i)}\sum^{n-1}_{h,p=1}\xi_{p}a^{p}_{h}\partial_{x_n}(c(dx_{h}))+\frac{i}{2(\xi_n-i)}\sum^{n-1}_{h=1}a^{n}_{h}\partial_{x_n}(c(dx_{h}))\nonumber\\
&+\frac{2i-\xi_n}{4(\xi_n-i)^2}h'(0)\sum^{n}_{h=1}\sum^{n-1}_{p=1}\xi_{p}a^{p}_{h}c(dx_{h})-\frac{1}{4(\xi_n-i)^2}h'(0)\sum^{n}_{h=1}a^{n}_{h}c(dx_{h}),\nonumber
\end{align}
where $\sum_{h=1}^{n-1}\partial_{x_n}(c(dx_h))=\sum_{h=1}^{n-1}\frac{1}{2}h'(0)c(dx_h).$\\
By calculation, we have
\begin{align}
\partial_{\xi_n}\left(\frac{ic[J(\xi)]}{|\xi|^{4}}\right)(x_0)
=\partial_{\xi_n}\left(\frac{i\sum^{n}_{i=1}\xi_{i}c[J(dx_i)]}{|\xi|^{4}}\right)(x_0)
=\partial_{\xi_n}\left(\frac{i\sum^{n}_{i,\beta=1}\xi_{i}a_{\beta}^ic(dx_{\beta})}{|\xi|^{4}}\right)(x_0),
\end{align}
\begin{align}
\partial_{\xi_n}\left(\frac{ic[J(\xi)]}{|\xi|^{4}}\right)(x_0)|_{|\xi'|=1}
&=i\sum^{n}_{\beta=1}\sum^{n-1}_{i=1}\xi_{i}a_{\beta}^ic(dx_{\beta})\partial_{\xi_n}\left(\frac{1}{(1+\xi_{n}^{2})^{2}}\right)+i\sum^{n}_{\beta=1}a_{\beta}^nc(dx_{\beta})\partial_{\xi_n}\left(\frac{\xi_{n}}{(1+\xi_{n}^{2})^{2}}\right)\\
&=-\frac{4 i \xi_n}{\left(\xi_n^2+1\right)^3}\sum^{n}_{\beta=1}\sum^{n-1}_{i=1}\xi_{i}a_{\beta}^ic(dx_{\beta})+\frac{i(1-3 \xi _n^2)}{\left(\xi _n^2+1\right){}^3}\sum^{n}_{\beta=1}a_{\beta}^nc(dx_{\beta})\nonumber
\end{align}
and
\begin{align}
\partial_{\xi_n}^2\left(\frac{ic[J(\xi)]}{|\xi|^{4}}\right)(x_0)|_{|\xi'|=1}
&=\frac{4 i \left(5 \xi _n^2-1\right)}{\left(\xi _n^2+1\right)^4}\sum^{n}_{\beta=1}\sum^{n-1}_{i=1}\xi_{i}a_{\beta}^ic(dx_{\beta})+\frac{12 i \xi _n \left(\xi _n^2-1\right)}{\left(\xi _n^2+1\right)^4}\sum^{n}_{\beta=1}a_{\beta}^nc(dx_{\beta}).
\end{align}
A trivial verification shows that
\begin{align}
&{\rm trace} [\partial_{x_n}\pi^+_{\xi_n}\sigma_{-1}({{D}_{J}}^{-1})\times
\partial_{\xi_n}^2\sigma_{-3}({{D}_{J}}^{-3})](x_0)|_{|\xi'|=1}\\
&=-\frac{2 i \left(5 \xi _n^2-1\right)}{\left(\xi _n-i\right)^5 \left(\xi _n+i\right)^4}\sum_{h=1}^{n}\sum_{p,i=1}^{n-1}\xi_{p}\xi_{i}a_{h}^{i}\partial_{x_n}(a_{h}^{p}){\rm tr}[\texttt{id}]
+\frac{6 \xi _n \left(\xi _n^2-1\right)}{\left(\xi _n-i\right)^5 \left(\xi _n+i\right)^4}\sum_{h=1}^{n}a_{h}^{n}\partial_{x_n}(a_{h}^{n}){\rm tr}[\texttt{id}]\nonumber\\
&-\frac{2-10 \xi _n^2}{\left(\xi _n-i\right)^5 \left(\xi _n+i\right)^4}\sum_{h=1}^{n}\sum_{i=1}^{n-1}\xi_{i}a_{h}^{i}\partial_{x_n}(a_{h}^{n}){\rm tr}[\texttt{id}]
-\frac{6 i \xi _n \left(\xi _n^2-1\right)}{\left(\xi _n-i\right)^5 \left(\xi _n+i\right)^4}\sum_{h=1}^{n}\sum_{p=1}^{n-1}\xi_{p}a_{h}^{n}\partial_{x_n}(a_{h}^{p}){\rm tr}[\texttt{id}]\nonumber\\
&-\frac{i \left(5 \xi _n^2-1\right)}{\left(\xi _n-i\right)^5 \left(\xi _n+i\right)^4}h'(0)\sum_{h,p,i=1}^{n-1}\xi_{p}\xi_{i}a_{h}^{p}a_{h}^{i}{\rm tr}[\texttt{id}]
-\frac{3 i \xi _n \left(\xi _n^2-1\right)}{\left(\xi _n-i\right)^5 \left(\xi _n+i\right)^4}h'(0)\sum_{h,p=1}^{n-1}\xi_{p}a_{h}^{p}a_{h}^{n}{\rm tr}[\texttt{id}]\nonumber\\
&-\frac{1-5 \xi _n^2}{\left(\xi _n-i\right)^5 \left(\xi _n+i\right)^4}h'(0)\sum_{h,i=1}^{n-1}\xi_{i}a_{h}^{n}a_{h}^{i}{\rm tr}[\texttt{id}]
+\frac{3 \xi _n \left(\xi _n^2-1\right)}{\left(\xi _n-i\right)^5 \left(\xi _n+i\right)^4}h'(0)\sum_{h=1}^{n-1}(a_{h}^{n})^2{\rm tr}[\texttt{id}]\nonumber\\
&+\frac{i \left(\xi _n-2 i\right) \left(5 \xi _n^2-1\right)}{\left(\xi _n-i\right)^6 \left(\xi _n+i\right)^4}h'(0)\sum_{h=1}^{n}\sum_{p,i=1}^{n-1}\xi_{p}\xi_{i}a_{h}^{p}a_{h}^{i}{\rm tr}[\texttt{id}]
+\frac{3 i \xi _n \left(\xi _n^2-1\right)}{\left(\xi _n-i\right)^6 \left(\xi _n+i\right)^4}h'(0)\sum_{h=1}^{n}(a_{h}^{n})^2{\rm tr}[\texttt{id}]\nonumber\\
&+\frac{i \left(5 \xi _n^2-1\right)}{\left(\xi _n-i\right)^6 \left(\xi _n+i\right)^4}h'(0)\sum_{h=1}^{n}\sum_{i=1}^{n-1}\xi_{i}a_{h}^{n}a_{h}^{i}{\rm tr}[\texttt{id}]
+\frac{3 i \xi _n \left(\xi _n-2 i\right) \left(\xi _n^2-1\right)}{\left(\xi _n-i\right)^6 \left(\xi _n+i\right)^4}h'(0)\sum_{h=1}^{n}\sum_{p=1}^{n-1}\xi_{p}a_{h}^{p}a_{h}^{n}{\rm tr}[\texttt{id}].\nonumber
\end{align}
Therefore
\begin{align}
\Phi_2
&=-\frac{1}{2}\int_{|\xi'|=1}\int^{+\infty}_{-\infty} {\rm trace}[\partial_{x_{n}}\pi^{+}_{\xi_{n}}\sigma_{-1}({{D}_{J}}^{-1})\times\partial^{2}_{\xi_{n}}\sigma_{-3}({{D}_{J}}^{-3})](x_0)d\xi_n\sigma(\xi')dx'\\
&=-\frac{1}{2}\int_{|\xi'|=1}\int^{+\infty}_{-\infty}-\frac{2 i \left(5 \xi _n^2-1\right)}{\left(\xi _n-i\right)^5 \left(\xi _n+i\right)^4}\sum_{h=1}^{n}\sum_{p,i=1}^{n-1}\xi_{p}\xi_{i}a_{h}^{i}\partial_{x_n}(a_{h}^{p}){\rm tr}[\texttt{id}]d\xi_n\sigma(\xi')dx'\nonumber\\
&-\frac{1}{2}\int_{|\xi'|=1}\int^{+\infty}_{-\infty}\frac{6 \xi _n \left(\xi _n^2-1\right)}{\left(\xi _n-i\right)^5 \left(\xi _n+i\right)^4}\sum_{h=1}^{n}a_{h}^{n}\partial_{x_n}(a_{h}^{n}){\rm tr}[\texttt{id}]d\xi_n\sigma(\xi')dx'\nonumber\\
&-\frac{1}{2}\int_{|\xi'|=1}\int^{+\infty}_{-\infty}-\frac{i \left(5 \xi _n^2-1\right)}{\left(\xi _n-i\right)^5 \left(\xi _n+i\right)^4}h'(0)\sum_{h,p,i=1}^{n-1}\xi_{p}\xi_{i}a_{h}^{p}a_{h}^{i}{\rm tr}[\texttt{id}]d\xi_n\sigma(\xi')dx'\nonumber\\
&-\frac{1}{2}\int_{|\xi'|=1}\int^{+\infty}_{-\infty}\frac{3 \xi _n \left(\xi _n^2-1\right)}{\left(\xi _n-i\right)^5 \left(\xi _n+i\right)^4}h'(0)\sum_{h=1}^{n-1}(a_{h}^{n})^2{\rm tr}[\texttt{id}]d\xi_n\sigma(\xi')dx'\nonumber
\end{align}
\begin{align}
&-\frac{1}{2}\int_{|\xi'|=1}\int^{+\infty}_{-\infty}\frac{i \left(\xi _n-2 i\right) \left(5 \xi _n^2-1\right)}{\left(\xi _n-i\right)^6 \left(\xi _n+i\right)^4}h'(0)\sum_{h=1}^{n}\sum_{p,i=1}^{n-1}\xi_{p}\xi_{i}a_{h}^{p}a_{h}^{i}{\rm tr}[\texttt{id}]d\xi_n\sigma(\xi')dx'\nonumber\\
&-\frac{1}{2}\int_{|\xi'|=1}\int^{+\infty}_{-\infty}\frac{3 i \xi _n \left(\xi _n^2-1\right)}{\left(\xi _n-i\right)^6 \left(\xi _n+i\right)^4}h'(0)\sum_{h=1}^{n}(a_{h}^{n})^2{\rm tr}[\texttt{id}]d\xi_n\sigma(\xi')dx'.\nonumber
\end{align}
As in the proof of  case (a)~(I), equation (3.53) gives
\begin{align}
\Phi_2
&=\sum_{h=1}^{n}\sum_{i=1}^{n-1}a_{h}^{i}\partial_{x_n}(a_{h}^{i}){\rm tr}[\texttt{id}]\Omega_4(-\frac{1}{2})\frac{8 \pi^2}{15}\frac{2\pi i}{4!}\Big[-\frac{2 i \left(5 \xi _n^2-1\right)}{\left(\xi _n+i\right)^4}\Big]^{(4)}\bigg|_{\xi_n=i}dx'\\
&+\sum_{h=1}^{n}a_{h}^{n}\partial_{x_n}(a_{h}^{n}){\rm tr}[\texttt{id}]\Omega_4(-\frac{1}{2})\frac{2\pi i}{4!}\Big[\frac{6 \xi _n \left(\xi _n^2-1\right)}{\left(\xi _n+i\right)^4}\Big]^{(4)}\bigg|_{\xi_n=i}dx'\nonumber\\
&+\sum_{h,i=1}^{n-1}(a_{h}^{i})^{2}{\rm tr}[\texttt{id}]\Omega_4 h'(0)(-\frac{1}{2})\frac{8 \pi^2}{15}\frac{2 \pi i}{4!}\Big[-\frac{i \left(5 \xi _n^2-1\right)}{\left(\xi _n+i\right)^4}\Big]^{(4)}\bigg|_{\xi_n=i}dx'\nonumber\\
&+\sum_{h=1}^{n-1}(a_{h}^{n})^2{\rm tr}[\texttt{id}]\Omega_4 h'(0)(-\frac{1}{2})\frac{2\pi i}{4!}\Big[\frac{3 \xi _n \left(\xi _n^2-1\right)}{\left(\xi _n+i\right)^4}\Big]^{(4)}\bigg|_{\xi_n=i}dx'\nonumber\\
&+\sum_{h=1}^{n}\sum_{i=1}^{n-1}(a_{h}^{i})^{2}{\rm tr}[\texttt{id}]\Omega_4 h'(0)(-\frac{1}{2})\frac{8 \pi^2}{15}\frac{2\pi i}{5!}\Big[\frac{i \left(\xi _n-2 i\right) \left(5 \xi _n^2-1\right)}{\left(\xi _n+i\right)^4}\Big]^{(5)}\bigg|_{\xi_n=i}dx'\nonumber\\
&+\sum_{h=1}^{n}(a_{h}^{n})^2{\rm tr}[\texttt{id}]\Omega_4h'(0)(-\frac{1}{2})\frac{2\pi i}{5!}\Big[\frac{3 i \xi _n \left(\xi _n^2-1\right)}{\left(\xi _n+i\right)^4}\Big]^{(5)}\bigg|_{\xi_n=i}dx'\nonumber\\
&=\sum_{h=1}^{n}\sum_{i=1}^{n-1}a_{h}^{i}\partial_{x_n}(a_{h}^{i}){\rm tr}[\texttt{id}]\Omega_4(\frac{\pi ^3}{24})dx'
+\sum_{h=1}^{n}a_{h}^{n}\partial_{x_n}(a_{h}^{n}){\rm tr}[\texttt{id}]\Omega_4(\frac{3 \pi }{64})dx'\nonumber\\
&+\sum_{h,i=1}^{n-1}(a_{h}^{i})^{2}{\rm tr}[\texttt{id}]\Omega_4 h'(0)(\frac{\pi ^3}{48})dx'
+\sum_{h=1}^{n-1}(a_{h}^{n})^2{\rm tr}[\texttt{id}]\Omega_4 h'(0)(\frac{3 \pi }{128})dx'\nonumber\\
&+\sum_{h=1}^{n}\sum_{i=1}^{n-1}(a_{h}^{i})^{2}{\rm tr}[\texttt{id}]\Omega_4 h'(0)(-\frac{7 \pi ^3}{120})dx'
+\sum_{h=1}^{n}(a_{h}^{n})^2{\rm tr}[\texttt{id}]\Omega_4h'(0)(-\frac{3 \pi}{64})dx'.\nonumber
\end{align}

\noindent  {\bf case (a)~(III)}~$r=-1,l=-3,|\alpha|=j=0,k=1$.\\

\noindent Using (3.18), we get
\begin{equation}
\Phi_3=-\frac{1}{2}\int_{|\xi'|=1}\int^{+\infty}_{-\infty}{\rm trace}[\partial_{\xi_{n}}\pi^{+}_{\xi_{n}}\sigma_{-1}({{D}_{J}}^{-1})
      \times\partial_{\xi_{n}}\partial_{x_{n}}\sigma_{-3}({{D}_{J}}^{-3})](x_0)d\xi_n\sigma(\xi')dx'.
\end{equation}\\
We check at once that
\begin{align}
\pi^+_{\xi_n}\partial_{\xi_n}\left(\frac{ic[J(\xi)]}{|\xi|^2}\right)(x_0)|_{|\xi'|=1}
&=-\frac{1}{2(\xi_{n}-i)^2}\sum^{n}_{\beta=1}\sum^{n-1}_{i=1}\xi_{i}a_{\beta}^ic(dx_{\beta})-\frac{i}{2(\xi_{n}-i)^2}\sum^{n}_{\beta=1}a_{\beta}^nc(dx_{\beta}).
\end{align}
Likewise,
\begin{align}
\partial_{x_n}\left(\frac{ic[J(\xi)]}{|\xi|^{4}}\right)(x_0)
&=\frac{i\partial_{x_n}(c[J(\xi)])(x_0)}{|\xi|^{4}}-\frac{ic[J(\xi)]\partial_{x_n}(|\xi|^{4})(x_0)}{|\xi|^{8}}\\
&=\frac{i\sum^{n}_{p,h=1}\xi_{p}\partial_{x_n}(a^{p}_{h})c(dx_{h})(x_0)}{|\xi|^{4}}+\frac{i\sum^{n}_{p,h=1}\xi_{p}a^{p}_{h}\partial_{x_n}(c(dx_{h}))(x_0)}{|\xi|^{4}}\nonumber\\
&-\frac{2ih'(0)|\xi^{'}|^{2}\sum^{n}_{p,h=1}\xi_{p}a^{p}_{h}c(dx_{h})(x_0)}{|\xi|^{6}}.\nonumber
\end{align}
Then we have
\begin{align}
\partial_{x_n}\left(\frac{ic[J(\xi)]}{|\xi|^{4}}\right)(x_0)|_{|\xi'|=1}
&=\frac{i\sum^{n}_{p,h=1}\xi_{p}\partial_{x_n}(a^{p}_{h})c(dx_{h})(x_0)}{(1+\xi_{n}^{2})^{2}}+\frac{i\sum^{n}_{p,h=1}\xi_{p}a^{p}_{h}\partial_{x_n}(c(dx_{h}))(x_0)}{(1+\xi_{n}^{2})^{2}}\\
&-\frac{2ih'(0)\sum^{n}_{p,h=1}\xi_{p}a^{p}_{h}c(dx_{h})(x_0)}{(1+\xi_{n}^{2})^{3}}.\nonumber
\end{align}
Thus
\begin{align}
&\partial_{\xi_n}\partial_{x_n}\left(\frac{ic[J(\xi)]}{|\xi|^4}\right)(x_0)|_{|\xi'|=1}\\
&=-\frac{4 i \xi _n}{\left(\xi _n^2+1\right)^3}\sum^{n}_{h=1}\sum^{n-1}_{p=1}\xi_{p}\partial_{x_n}(a^{p}_{h})c(dx_{h})
+\frac{i(1-3 \xi _n^2)}{\left(\xi _n^2+1\right)^3}\sum^{n}_{h=1}\partial_{x_n}(a^{n}_{h})c(dx_{h})\nonumber\\
&-\frac{4 i \xi _n}{\left(\xi _n^2+1\right)^3}\sum^{n-1}_{h,p=1}\xi_{p}a^{p}_{h}\partial_{x_n}(c(dx_{h}))
+\frac{i(1-3 \xi _n^2)}{\left(\xi _n^2+1\right)^3}\sum^{n-1}_{h=1}a^{n}_{h}\partial_{x_n}(c(dx_{h}))\nonumber\\
&+\frac{12 i \xi _n}{\left(\xi _n^2+1\right)^4}h'(0)\sum^{n}_{h=1}\sum^{n-1}_{p=1}\xi_{p}a^{p}_{h}c(dx_{h})
-\frac{2i(1-5 \xi _n^2)}{\left(\xi _n^2+1\right)^4}h'(0)\sum^{n}_{h=1}a^{n}_{h}c(dx_{h}).\nonumber
\end{align}
We see at once that
\begin{align}
&{\rm trace} [\partial_{\xi_n}\pi^+_{\xi_n}\sigma_{-1}({{D}_{J}}^{-1})\times \partial_{\xi_n}\partial_{x_n}\sigma_{-3}({{D}_{J}}^{-3})](x_0)|_{|\xi'|=1}\\
&=-\frac{2 i \xi _n}{\left(\xi _n-i\right)^5 \left(\xi _n+i\right)^3}\sum_{h=1}^{n}\sum_{i,p=1}^{n-1}\xi_{i}\xi_{p}a_{h}^{i}\partial_{x_n}(a_{h}^{p}){\rm tr}[\texttt{id}]
-\frac{i \left(3 \xi _n^2-1\right)}{4 \left(\xi _n-i\right)^5 \left(\xi _n+i\right)^3}h'(0)\sum_{i,h=1}^{n-1}\xi_{i}a_{h}^{i}a_{h}^{n}{\rm tr}[\texttt{id}]\nonumber\\
&-\frac{i \xi _n}{\left(\xi _n-i\right)^5 \left(\xi _n+i\right)^3}h'(0)\sum_{i,h,p=1}^{n-1}\xi_{i}\xi_{p}a_{h}^{i}a_{h}^{p}{\rm tr}[\texttt{id}]
-\frac{i \left(3 \xi _n^2-1\right)}{2 \left(\xi _n-i\right)^5 \left(\xi _n+i\right)^3}\sum_{h=1}^{n}\sum_{i=1}^{n-1}\xi_{i}a_{h}^{i}\partial_{x_n}(a_{h}^{n}){\rm tr}[\texttt{id}]\nonumber\\
&+\frac{6 i \xi _n}{\left(\xi _n-i\right)^6 \left(\xi _n+i\right)^4}h'(0)\sum_{h=1}^{n}\sum_{i,p=1}^{n-1}\xi_{i}\xi_{p}a_{h}^{i}a_{h}^{p}{\rm tr}[\texttt{id}]
-\frac{1-3 \xi _n^2}{2 \left(\xi _n-i\right)^5 \left(\xi _n+i\right)^3}\sum_{h=1}^{n}a_{h}^{n}\partial_{x_n}(a_{h}^{n}){\rm tr}[\texttt{id}]\nonumber\\
&+\frac{2 \xi _n}{\left(\xi _n-i\right)^5 \left(\xi _n+i\right)^3}\sum_{h=1}^{n}\sum_{p=1}^{n-1}\xi_{p}a_{h}^{n}\partial_{x_n}(a_{h}^{p}){\rm tr}[\texttt{id}]
-\frac{i \left(1-5 \xi _n^2\right)}{\left(\xi _n-i\right)^6 \left(\xi _n+i\right)^4}h'(0)\sum_{h=1}^{n}\sum_{i=1}^{n-1}\xi_{i}a_{h}^{i}a_{h}^{n}{\rm tr}[\texttt{id}]\nonumber\\
&+\frac{\xi _n}{\left(\xi _n-i\right)^5 \left(\xi _n+i\right)^3}h'(0)\sum_{h,p=1}^{n-1}\xi_{p}a_{h}^{n}a_{h}^{p}{\rm tr}[\texttt{id}]
-\frac{1-3 \xi _n^2}{4 \left(\xi _n-i\right)^5 \left(\xi _n+i\right)^3}h'(0)\sum_{h=1}^{n-1}(a_{h}^{n})^2{\rm tr}[\texttt{id}]\nonumber\\
&-\frac{6 \xi _n}{\left(\xi _n-i\right)^6 \left(\xi _n+i\right)^4}h'(0)\sum_{h=1}^{n}\sum_{p=1}^{n-1}\xi_{p}a_{h}^{n}a_{h}^{p}{\rm tr}[\texttt{id}]
-\frac{5 \xi _n^2-1}{\left(\xi _n-i\right)^6 \left(\xi _n+i\right)^4}h'(0)\sum_{h=1}^{n}(a_{h}^{n})^2{\rm tr}[\texttt{id}].\nonumber
\end{align}
From this, we have
\begin{align}
\Phi_3
&=-\frac{1}{2}\int_{|\xi'|=1}\int^{+\infty}_{-\infty}
{\rm trace} [\partial_{\xi_n}\pi^+_{\xi_n}\sigma_{-1}({{D}_{J}}^{-1})\times
\partial_{\xi_n}\partial_{x_n}\sigma_{-3}({{D}_{J}}^{-3})](x_0)d\xi_n\sigma(\xi')dx'\\
&=-\frac{1}{2}\int_{|\xi'|=1}\int^{+\infty}_{-\infty}-\frac{2 i \xi _n}{\left(\xi _n-i\right)^5 \left(\xi _n+i\right)^3}\sum_{h=1}^{n}\sum_{i,p=1}^{n-1}\xi_{i}\xi_{p}a_{h}^{i}\partial_{x_n}(a_{h}^{p}){\rm tr}[\texttt{id}]d\xi_n\sigma(\xi')dx'\nonumber\\
&-\frac{1}{2}\int_{|\xi'|=1}\int^{+\infty}_{-\infty}-\frac{i \xi _n}{\left(\xi _n-i\right)^5 \left(\xi _n+i\right)^3}h'(0)\sum_{i,h,p=1}^{n-1}\xi_{i}\xi_{p}a_{h}^{i}a_{h}^{p}{\rm tr}[\texttt{id}]d\xi_n\sigma(\xi')dx'\nonumber\\
&-\frac{1}{2}\int_{|\xi'|=1}\int^{+\infty}_{-\infty}\frac{6 i \xi _n}{\left(\xi _n-i\right)^6 \left(\xi _n+i\right)^4}h'(0)\sum_{h=1}^{n}\sum_{i,p=1}^{n-1}\xi_{i}\xi_{p}a_{h}^{i}a_{h}^{p}{\rm tr}[\texttt{id}]d\xi_n\sigma(\xi')dx'\nonumber\\
&-\frac{1}{2}\int_{|\xi'|=1}\int^{+\infty}_{-\infty}-\frac{1-3 \xi _n^2}{2 \left(\xi _n-i\right)^5 \left(\xi _n+i\right)^3}\sum_{h=1}^{n}a_{h}^{n}\partial_{x_n}(a_{h}^{n}){\rm tr}[\texttt{id}]d\xi_n\sigma(\xi')dx'\nonumber
\end{align}
\begin{align}
&-\frac{1}{2}\int_{|\xi'|=1}\int^{+\infty}_{-\infty}-\frac{1-3 \xi _n^2}{4 \left(\xi _n-i\right)^5 \left(\xi _n+i\right)^3}h'(0)\sum_{h=1}^{n-1}(a_{h}^{n})^2{\rm tr}[\texttt{id}]d\xi_n\sigma(\xi')dx'\nonumber\\
&-\frac{1}{2}\int_{|\xi'|=1}\int^{+\infty}_{-\infty}-\frac{5 \xi _n^2-1}{\left(\xi _n-i\right)^6 \left(\xi _n+i\right)^4}h'(0)\sum_{h=1}^{n}(a_{h}^{n})^2{\rm tr}[\texttt{id}]d\xi_n\sigma(\xi')dx'.\nonumber
\end{align}
It is immediate that
\begin{align}
\Phi_3
&=\sum_{h=1}^{n}\sum_{i=1}^{n-1}a_{h}^{i}\partial_{x_n}(a_{h}^{i}){\rm tr}[\texttt{id}]\Omega_4(-\frac{1}{2})\frac{8 \pi^{2}}{15}\frac{2\pi i}{4!}\Big[-\frac{2 i \xi _n}{\left(\xi _n+i\right)^3}\Big]^{(4)}\bigg|_{\xi_n=i}dx'\\
&+\sum_{i,h=1}^{n-1}(a_{h}^{i})^{2}{\rm tr}[\texttt{id}]\Omega_4h'(0)(-\frac{1}{2})\frac{8 \pi^{2}}{15}\frac{2\pi i}{4!}\Big[-\frac{i \xi _n}{\left(\xi _n+i\right)^3}\Big]^{(4)}\bigg|_{\xi_n=i}dx'\nonumber\\
&+\sum_{h=1}^{n}\sum_{i=1}^{n-1}(a_{h}^{i})^{2}{\rm tr}[\texttt{id}]\Omega_4h'(0)(-\frac{1}{2})\frac{8 \pi^{2}}{15}\frac{2\pi i}{5!}\Big[\frac{6 i \xi _n}{\left(\xi _n+i\right)^4}\Big]^{(5)}\bigg|_{\xi_n=i}dx'\nonumber\\
&+\sum_{h=1}^{n}a_{h}^{n}\partial_{x_n}(a_{h}^{n}){\rm tr}[\texttt{id}]\Omega_4(-\frac{1}{2})\frac{2\pi i}{4!}\Big[-\frac{1-3 \xi _n^2}{2 \left(\xi _n+i\right)^3}\Big]^{(4)}\bigg|_{\xi_n=i}dx'\nonumber\\
&+\sum_{h=1}^{n-1}(a_{h}^{n})^2{\rm tr}[\texttt{id}]\Omega_4h'(0)(-\frac{1}{2})\frac{2\pi i}{4!}\Big[-\frac{1-3 \xi _n^2}{4 \left(\xi _n+i\right)^3}\Big]^{(4)}\bigg|_{\xi_n=i}dx'\nonumber\\
&+\sum_{h=1}^{n}(a_{h}^{n})^2{\rm tr}[\texttt{id}]\Omega_4h'(0)(-\frac{1}{2})\frac{2\pi i}{5!}\Big[-\frac{5 \xi _n^2-1}{\left(\xi _n+i\right)^4}\Big]^{(5)}\bigg|_{\xi_n=i}dx'\nonumber\\
&=\sum_{h=1}^{n}\sum_{i=1}^{n-1}a_{h}^{i}\partial_{x_n}(a_{h}^{i}){\rm tr}[\texttt{id}]\Omega_4(-\frac{\pi ^3}{24})dx'
+\sum_{i,h=1}^{n-1}(a_{h}^{i})^{2}{\rm tr}[\texttt{id}]\Omega_4h'(0)(-\frac{\pi ^3}{48} )dx'\nonumber\\
&+\sum_{h=1}^{n}\sum_{i=1}^{n-1}(a_{h}^{i})^{2}{\rm tr}[\texttt{id}]\Omega_4h'(0)(\frac{7 \pi ^3}{80})dx'
+\sum_{h=1}^{n}a_{h}^{n}\partial_{x_n}(a_{h}^{n}){\rm tr}[\texttt{id}]\Omega_4(-\frac{3 \pi}{64} )dx'\nonumber\\
&+\sum_{h=1}^{n-1}(a_{h}^{n})^2{\rm tr}[\texttt{id}]\Omega_4h'(0)(-\frac{3 \pi}{128})dx'
+\sum_{h=1}^{n}(a_{h}^{n})^2{\rm tr}[\texttt{id}]\Omega_4h'(0)(\frac{9 \pi }{128})dx'.\nonumber
\end{align}

We have the facts that
\begin{align}
&a^{p}_{l}a^{j}_{p}=\delta_l^j, a^{p}_{l}=a^{l}_{p};
\end{align}
\begin{align}
&\sum_{h=1}^{n}\sum_{i=1}^{n-1}a_{h}^{i}\partial_{x_i}(a_{h}^{n})=\sum_{h,i=1}^{n}a_{h}^{i}\partial_{x_i}(a_{h}^{n})
\end{align}
and
\begin{align}
\sum_{h=1}^{n}\sum_{i=1}^{n-1}a_{h}^{i}\partial_{x_i}(a_{h}^{n})=-\sum_{h=1}^{n}\sum_{i=1}^{n-1}a_{h}^{n}\partial_{x_i}(a_{h}^{i}).
\end{align}

In combination with the calculation,
\begin{align}
\Phi_1+\Phi_2+\Phi_3&=
\sum_{l,j=1}^{n}a_{l}^{j}\partial_{x_j}(a_{l}^{n}){\rm tr}[\texttt{id}]\Omega_4(-\frac{\pi}{16}+\frac{\pi^3}{6})dx'\\
&+\sum_{l=1}^{n}\sum_{i=1}^{n-1}(a_{l}^{i})^2{\rm tr}[\texttt{id}]\Omega_4h'(0)(\frac{7 \pi ^3}{240})dx'\nonumber\\
&+\sum_{l=1}^{n}(a_{l}^{n})^2{\rm tr}[\texttt{id}]\Omega_4h'(0)(\frac{3 \pi }{128})dx'.\nonumber
\end{align}

\noindent  {\bf case (b)}~$r=-1,l=-4,|\alpha|=j=k=0$.\\

\noindent It is easily seen that
\begin{align}
\Phi_4&=-i\int_{|\xi'|=1}\int^{+\infty}_{-\infty}{\rm trace}[\pi^{+}_{\xi_{n}}\sigma_{-1}({{D}_{J}}^{-1})
      \times\partial_{\xi_{n}}\sigma_{-4}({{D}_{J}}^{-3})](x_0)d\xi_n\sigma(\xi')dx'\\
&=i\int_{|\xi'|=1}\int^{+\infty}_{-\infty}{\rm trace} [\partial_{\xi_n}\pi^+_{\xi_n}\sigma_{-1}({{D}_{J}}^{-1})\times
\sigma_{-4}({{D}_{J}}^{-3})](x_0)d\xi_n\sigma(\xi')dx'.\nonumber
\end{align}
In the normal coordinate, $g^{ij}(x_{0})=\delta^{j}_{i}$ and $\partial_{x_{j}}(g^{\alpha\beta})(x_{0})=0$, if $j<n$; $\partial_{x_{j}}(g^{\alpha\beta})(x_{0})=h'(0)\delta^{\alpha}_{\beta}$, if $j=n$.
So by \cite{Wa3}, when $k<n$, we have $\Gamma^{n}(x_{0})=\frac{5}{2}h'(0)$ and $\Gamma^{k}(x_{0})=0.$
We thus get
\begin{align}
&\sigma_{-4}({{D}_{J}}^{-3})(x_{0})|_{|\xi'|=1}\\
&=\frac{1}{(1+\xi_{n}^{2})^{3}}h'(0)\sum_{\eta, \Gamma, \Omega, \Lambda, \Pi=1}^{n}\xi_{\Gamma}\xi_{\Omega}a_{\Lambda}^{\Gamma}a_{\eta}^{n}a_{\Pi}^{\Omega}c(dx_{\Lambda})c(dx_{\eta})c(dx_{\Pi})\nonumber\\
&-\frac{1}{(1+\xi_{n}^{2})^{3}}h'(0)\sum_{\chi,\tau=1}^{n}\sum_{\gamma=1}^{n-1}\xi_{\gamma}\xi_{\chi}a_{\tau}^{\chi}c(dx_{\gamma})c(dx_{n})c(dx_{\tau})\nonumber\\
&+\frac{5}{(1+\xi_{n}^{2})^{3}}h'(0)\sum_{\rho,\theta=1}^{n}\xi_{n}\xi_{\rho}a_{\theta}^{\rho}c(dx_{\theta})\nonumber\\
&+\frac{2}{(1+\xi_{n}^{2})^{3}}\sum_{\alpha,\beta,\lambda,\omega=1}^{n}\xi_{\lambda}a_{\alpha}^{\beta}a_{\omega}^{\lambda}c(dx_{\beta})c[(\nabla^{L}_{e_{\alpha}}J)(\xi^{*})]c(dx_{\omega})\nonumber\\
&-\frac{1}{4(1+\xi_{n}^{2})^{3}}h'(0)\sum_{\mu, \Phi, b, \Psi, c=1}^{n}\sum_{\nu=1}^{n-1}\xi_{\Phi}\xi_{b}a_{\nu}^{\mu}a_{\Psi}^{\Phi}a_{c}^{b}c(dx_{\Psi})c(dx_{\mu})c(dx_{n})c(dx_{\nu})c(dx_{c})\nonumber\\
&+\frac{1}{(1+\xi_{n}^{2})^{3}}\sum_{j,p,h,\delta,\varepsilon,q=1}^{n}\xi_{p}\xi_{\delta}a_{\varepsilon}^{\delta}a_{q}^{j}\partial_{x_{j}}(a_{h}^{p})c(dx_{\varepsilon})c(dx_{q})c(dx_{h})\nonumber\\
&-\frac{2}{(1+\xi_{n}^{2})^{3}}\sum_{j,p,h=1}^{n}\xi_{j}\xi_{p}\partial_{x_{j}}(a_{h}^{p})c(dx_{h})\nonumber\\
&+\frac{1}{(1+\xi_{n}^{2})^{3}}\sum_{p,\kappa,o,e=1}^{n}\sum_{h=1}^{n-1}\xi_{p}\xi_{\kappa}a_{h}^{p}a_{o}^{\kappa}a_{e}^{n}c(dx_{o})c(dx_{e})\partial_{x_{n}}(c(dx_{h}))\nonumber\\
&-\frac{2}{(1+\xi_{n}^{2})^{3}}\sum_{p=1}^{n}\sum_{h=1}^{n-1}\xi_{n}\xi_{p}a_{h}^{p}\partial_{x_{n}}(c(dx_{h}))\nonumber\\
&-\frac{2}{(1+\xi_{n}^{2})^{3}}h'(0)\sum_{d, f, e, m, g=1}^{n}\xi_{d}\xi_{f}a_{e}^{d}a_{m}^{n}a_{g}^{f}c(dx_{e})c(dx_{m})c(dx_{g})\nonumber\\
&+\frac{4}{(1+\xi_{n}^{2})^{4}}h'(0)\sum_{\psi,\varphi=1}^{n}\xi_{n}\xi_{\psi}a_{\varphi}^{\psi}c(dx_{\varphi}).\nonumber
\end{align}
With the notation
\begin{align}
B_1(x_0)&=-\frac{1}{2(\xi_{n}-i)^2}\sum^{n}_{l=1}\sum^{n-1}_{i=1}\xi_{i}a_{l}^ic(dx_{l});\\
B_2(x_0)&=-\frac{i}{2(\xi_{n}-i)^2}\sum^{n}_{l=1}a_{l}^nc(dx_{l}),
\end{align}
we have
\begin{align}
i\int_{|\xi'|=1}\int^{+\infty}_{-\infty}{\rm trace} [B_1(x_0)\times
\sigma_{-4}({{D}_{J}}^{-3})](x_0)d\xi_n\sigma(\xi')dx'=0
\end{align}
and
\begin{align}
&i\int_{|\xi'|=1}\int^{+\infty}_{-\infty}{\rm trace} [B_2(x_0)\times
\sigma_{-4}({{D}_{J}}^{-3})](x_0)d\xi_n\sigma(\xi')dx'\\
&=i\int_{|\xi'|=1}\int^{+\infty}_{-\infty}-\frac{i}{2 \left(\xi _n-i\right)^5 \left(\xi _n+i\right)^3}h'(0)
\Big(\sum_{l, \eta, \Gamma, \Omega=1}^{n}\xi_{\Gamma}\xi_{\Omega} a_{l}^{n}a_{\eta}^{\Gamma}a_{\eta}^{n}a_{l}^{\Omega}{\rm tr}[\texttt{id}]\nonumber\\
&-\sum_{l, \Gamma, \Omega, \Lambda=1}^{n}\xi_{\Gamma}\xi_{\Omega} a_{l}^{n}a_{\Lambda}^{\Gamma}a_{l}^{n}a_{\Lambda}^{\Omega}{\rm tr}[\texttt{id}]+\sum_{l, \eta, \Gamma, \Omega=1}^{n}\xi_{\Gamma}\xi_{\Omega} a_{l}^{n}a_{l}^{\Gamma}a_{\eta}^{n}a_{\eta}^{\Omega}{\rm tr}[\texttt{id}]\Big)d\xi_n\sigma(\xi')dx'\nonumber\\
&+i\int_{|\xi'|=1}\int^{+\infty}_{-\infty}\frac{i}{2 \left(\xi _n-i\right)^5 \left(\xi _n+i\right)^3}h'(0)\Big(-\sum_{\chi=1}^{n}\sum_{\gamma=1}^{n-1}\xi_{\gamma}\xi_{\chi}a_{n}^{n}a_{\gamma}^{\chi}{\rm tr}[\texttt{id}]\nonumber\\
&+\sum_{\chi=1}^{n}\sum_{l=1}^{n-1}\xi_{l}\xi_{\chi}a_{l}^{n}a_{n}^{\chi}{\rm tr}[\texttt{id}]\Big)d\xi_n\sigma(\xi')dx'\nonumber\\
&+i\int_{|\xi'|=1}\int^{+\infty}_{-\infty}\frac{5i}{2 \left(\xi _n-i\right)^5 \left(\xi _n+i\right)^3}h'(0)\sum_{l,\rho=1}^{n}\xi_{n}\xi_{\rho}a_{l}^{n}a_{l}^{\rho}{\rm tr}[\texttt{id}]d\xi_n\sigma(\xi')dx'\nonumber\\
&+i\int_{|\xi'|=1}\int^{+\infty}_{-\infty}-\frac{i}{\left(\xi _n-i\right)^5 \left(\xi _n+i\right)^3}\Big(\sum_{l,\alpha,\beta,\lambda=1}^{n}\xi_{\lambda}a_{l}^{n}a_{\alpha}^{\beta}a_{l}^{\lambda}g^{M}(dx_{\beta}, (\nabla^{L}_{e_{\alpha}}J)(\xi^{*})){\rm tr}[\texttt{id}]\nonumber\\
&-\sum_{l,\alpha,\beta,\lambda=1}^{n}\xi_{\lambda}a_{l}^{n}a_{\alpha}^{\beta}a_{\beta}^{\lambda}g^{M}(dx_{l}, (\nabla^{L}_{e_{\alpha}}J)(\xi^{*})){\rm tr}[\texttt{id}]\nonumber\\
&+\sum_{l,\alpha,\lambda,\omega=1}^{n}\xi_{\lambda}a_{l}^{n}a_{\alpha}^{l}a_{\omega}^{\lambda}g^{M}(dx_{\omega}, (\nabla^{L}_{e_{\alpha}}J)(\xi^{*})){\rm tr}[\texttt{id}]\Big)d\xi_n\sigma(\xi')dx'\nonumber\\
&+i\int_{|\xi'|=1}\int^{+\infty}_{-\infty}\frac{i}{8 \left(\xi _n-i\right)^5 \left(\xi _n+i\right)^3}h'(0)\Big(-\sum_{l, \Phi, b=1}^{n}\sum_{\nu=1}^{n-1}\xi_{\Phi}\xi_{b}a_{l}^{n}a_{\nu}^{n}a_{\nu}^{\Phi}a_{l}^{b}{\rm tr}[\texttt{id}]\nonumber\\
&+\sum_{l, \Phi, b=1}^{n}\sum_{\nu=1}^{n-1}\xi_{\Phi}\xi_{b}a_{l}^{n}a_{\nu}^{\nu}a_{n}^{\Phi}a_{l}^{b}{\rm tr}[\texttt{id}]+\sum_{\Phi, b, \Psi=1}^{n}\sum_{\nu=1}^{n-1}\xi_{\Phi}\xi_{b}a_{\nu}^{n}a_{\nu}^{n}a_{\Psi}^{\Phi}a_{\Psi}^{b}{\rm tr}[\texttt{id}]\nonumber
\end{align}
\begin{align}
&-\sum_{\mu, \Phi, b=1}^{n}\sum_{\nu=1}^{n-1}\xi_{\Phi}\xi_{b}a_{\nu}^{n}a_{\nu}^{\mu}a_{n}^{\Phi}a_{\mu}^{b}{\rm tr}[\texttt{id}]+\sum_{\mu, \Phi, b=1}^{n}\sum_{\nu=1}^{n-1}\xi_{\Phi}\xi_{b}a_{\nu}^{n}a_{\nu}^{\mu}a_{\mu}^{\Phi}a_{n}^{b}{\rm tr}[\texttt{id}]\nonumber\\
&-\sum_{\Phi, b, \Psi=1}^{n}\sum_{\nu=1}^{n-1}\xi_{\Phi}\xi_{b}a_{n}^{n}a_{\nu}^{\nu}a_{\Psi}^{\Phi}a_{\Psi}^{b}{\rm tr}[\texttt{id}]+\sum_{\mu, \Phi, b=1}^{n}\sum_{\nu=1}^{n-1}\xi_{\Phi}\xi_{b}a_{n}^{n}a_{\nu}^{\mu}a_{\nu}^{\Phi}a_{\mu}^{b}{\rm tr}[\texttt{id}]\nonumber\\
&-\sum_{\mu, \Phi, b=1}^{n}\sum_{\nu=1}^{n-1}\xi_{\Phi}\xi_{b}a_{n}^{n}a_{\nu}^{\mu}a_{\mu}^{\Phi}a_{\nu}^{b}{\rm tr}[\texttt{id}]-\sum_{l, \Phi, b=1}^{n}\sum_{\nu=1}^{n-1}\xi_{\Phi}\xi_{b}a_{l}^{n}a_{\nu}^{l}a_{\nu}^{\Phi}a_{n}^{b}{\rm tr}[\texttt{id}]\nonumber\\
&+\sum_{l, \Phi, b=1}^{n}\sum_{\nu=1}^{n-1}\xi_{\Phi}\xi_{b}a_{l}^{n}a_{\nu}^{l}a_{n}^{\Phi}a_{\nu}^{b}{\rm tr}[\texttt{id}]+\sum_{l, \Phi, b=1}^{n}\sum_{\nu=1}^{n-1}\xi_{\Phi}\xi_{b}a_{l}^{n}a_{\nu}^{\nu}a_{l}^{\Phi}a_{n}^{b}{\rm tr}[\texttt{id}]\nonumber\\
&-\sum_{l, \Phi, b=1}^{n}\sum_{\nu=1}^{n-1}\xi_{\Phi}\xi_{b}a_{l}^{n}a_{\nu}^{n}a_{l}^{\Phi}a_{\nu}^{b}{\rm tr}[\texttt{id}]\Big)d\xi_n\sigma(\xi')dx'\nonumber\\
&+i\int_{|\xi'|=1}\int^{+\infty}_{-\infty}-\frac{i}{2 \left(\xi _n-i\right)^5 \left(\xi _n+i\right)^3}\Big(\sum_{l,j,p,\delta,q=1}^{n}\xi_{p}\xi_{\delta}a_{l}^{n}a_{q}^{\delta}a_{q}^{j}\partial_{x_{j}}(a_{l}^{p}){\rm tr}[\texttt{id}]\nonumber\\
&-\sum_{l,j,p,h,\delta=1}^{n}\xi_{p}\xi_{\delta}a_{l}^{n}a_{h}^{\delta}a_{l}^{j}\partial_{x_{j}}(a_{h}^{p}){\rm tr}[\texttt{id}]+\sum_{l,j,p,\delta,q=1}^{n}\xi_{p}\xi_{\delta}a_{l}^{n}a_{l}^{\delta}a_{q}^{j}\partial_{x_{j}}(a_{q}^{p}){\rm tr}[\texttt{id}]\Big)d\xi_n\sigma(\xi')dx'\nonumber\\
&+i\int_{|\xi'|=1}\int^{+\infty}_{-\infty}-\frac{i}{\left(\xi _n-i\right)^5 \left(\xi _n+i\right)^3}\sum_{l,j,p=1}^{n}\xi_{j}\xi_{p}a_{l}^{n}\partial_{x_{j}}(a_{l}^{p}){\rm tr}[\texttt{id}]d\xi_n\sigma(\xi')dx'\nonumber\\
&+i\int_{|\xi'|=1}\int^{+\infty}_{-\infty}-\frac{i}{4 \left(\xi _n-i\right)^5 \left(\xi _n+i\right)^3}h'(0)\Big(\sum_{p,\kappa,o=1}^{n}\sum_{l=1}^{n-1}\xi_{p}\xi_{\kappa}a_{l}^{n}a_{l}^{p}a_{o}^{\kappa}a_{o}^{n}{\rm tr}[\texttt{id}]\nonumber\\
&-\sum_{l,p,\kappa=1}^{n}\sum_{h=1}^{n-1}\xi_{p}\xi_{\kappa}a_{l}^{n}a_{h}^{p}a_{h}^{\kappa}a_{l}^{n}{\rm tr}[\texttt{id}]+\sum_{l,p,\kappa=1}^{n}\sum_{h=1}^{n-1}\xi_{p}\xi_{\kappa}a_{l}^{n}a_{h}^{p}a_{l}^{\kappa}a_{h}^{n}{\rm tr}[\texttt{id}]\Big)d\xi_n\sigma(\xi')dx'\nonumber\\
&+i\int_{|\xi'|=1}\int^{+\infty}_{-\infty}-\frac{i}{2\left(\xi _n-i\right)^5 \left(\xi _n+i\right)^3}h'(0)\sum_{p=1}^{n}\sum_{l=1}^{n-1}\xi_{n}\xi_{p}a_{l}^{n}a_{l}^{p}{\rm tr}[\texttt{id}]d\xi_n\sigma(\xi')dx'\nonumber\\
&+i\int_{|\xi'|=1}\int^{+\infty}_{-\infty}\frac{i}{\left(\xi _n-i\right)^6 \left(\xi _n+i\right)^4}h'(0)\Big(\sum_{l, d, f, e=1}^{n}\xi_{d}\xi_{f}a_{l}^{n}a_{e}^{d}a_{e}^{n}a_{l}^{f}{\rm tr}[\texttt{id}]\nonumber\\
&-\sum_{l, d, f, e=1}^{n}\xi_{d}\xi_{f}a_{l}^{n}a_{e}^{d}a_{l}^{n}a_{e}^{f}{\rm tr}[\texttt{id}]+\sum_{l, d, f, m=1}^{n}\xi_{d}\xi_{f}a_{l}^{n}a_{l}^{d}a_{m}^{n}a_{m}^{f}{\rm tr}[\texttt{id}]\Big)d\xi_n\sigma(\xi')dx'\nonumber\\
&+i\int_{|\xi'|=1}\int^{+\infty}_{-\infty}\frac{2i}{\left(\xi _n-i\right)^6 \left(\xi _n+i\right)^4}h'(0)\sum_{l,\psi=1}^{n}\xi_{n}\xi_{\psi}a_{l}^{n}a_{l}^{\psi}{\rm tr}[\texttt{id}]d\xi_n\sigma(\xi')dx'.\nonumber
\end{align}
Summarizing, we have
\begin{align}
\Phi_4
&=\Big(\sum_{l=1}^{n}g^{M}(J(dx_{l}), (\nabla^{L}_{e_{l}}J)e_{n}){\rm tr}[\texttt{id}]
-\sum_{l=1}^{n}g^{M}(J(dx_{n}), (\nabla^{L}_{e_{l}}J)e_{l}){\rm tr}[\texttt{id}]\nonumber\\
&+\sum_{l=1}^{n}g^{M}(J(dx_{l}), (\nabla^{L}_{e_{n}}J)e_{l}){\rm tr}[\texttt{id}]\Big)
\Omega_4(-\frac{\pi^{3}}{8})dx'\nonumber\\
&+\sum_{l,j=1}^{n}a_{l}^{j}\partial_{x_{j}}(a_{l}^{n}){\rm tr}[\texttt{id}]\Omega_4(\frac{\pi^{3}}{8})dx'\nonumber\\
&+\sum_{l,\beta=1}^{n}\sum_{i=1}^{n-1}(a_{l}^{n})^{2}(a_{\beta}^{i})^{2}{\rm tr}[\texttt{id}]\Omega_4h'(0)(\frac{\pi^{3}}{32})dx'\nonumber\\
&+\sum_{l=1}^{n}(a_{l}^{n})^{2}{\rm tr}[\texttt{id}]\Omega_4h'(0)(\frac{3\pi^{3}}{5})dx'\nonumber\\
&+\sum_{i=1}^{n-1}(a_{i}^{n})^{2}{\rm tr}[\texttt{id}]\Omega_4h'(0)(-\frac{5\pi^{3}}{64})dx'\nonumber\\
&+\sum_{i=1}^{n-1}a_{n}^{n}a_{i}^{i}{\rm tr}[\texttt{id}]\Omega_4h'(0)(-\frac{3\pi^{3}}{64})dx'.\nonumber
\end{align}

\noindent {\bf  case (c)}~$r=-2,l=-3,|\alpha|=j=k=0$.\\

\noindent We calculate
\begin{equation}
\Phi_5=-i\int_{|\xi'|=1}\int^{+\infty}_{-\infty}{\rm trace} [\pi^{+}_{\xi_{n}}\sigma_{-2}({{D}_{J}}^{-1})
      \times\partial_{\xi_{n}}\sigma_{-3}({{D}_{J}}^{-3})](x_0)d\xi_n\sigma(\xi')dx'.
\end{equation}\\
We follow the notation of \cite{LW2}.
\begin{align}
A_1(x_0)&=\frac{c[J(\xi)]\sigma_{0}({D}_{J})(x_0)c[J(\xi)]}{(1+\xi_n^2)^2};\\
A_2(x_0)&=\frac{c[J(\xi)]}{(1+\xi_n^2)^2}\Big[\sum_{j,p,h=1}^{n}\xi_p\partial_{x_j}(a_{h}^{p})c[J(dx_j)]c(dx_h)+\sum_{p=1}^{n}\sum_{h=1}^{n-1}\xi_pa_{h}^{p}c[J(dx_n)]\partial_{x_n}(c(dx_h))\Big];\\
A_3(x_0)&=\frac{c[J(\xi)]}{(1+\xi_n^2)^3}c[J(dx_n)]c[J(\xi)],
\end{align}
means that
\begin{align}
&\pi^+_{\xi_n}\sigma_{-2}({{D}_{J}}^{-1})(x_0)|_{|\xi'|=1}=\pi^+_{\xi_n}(A_1(x_0))+\pi^+_{\xi_n}(A_2(x_0))-h'(0)\pi^+_{\xi_n}(A_3(x_0)).
\end{align}
Computations show that
\begin{align}
\pi^+_{\xi_n}(A_1(x_0))
&=\frac{i\xi_n}{16(\xi_n-i)^2}h'(0)\sum_{l,\gamma,\mu=1}^{n}\sum_{\nu=1}^{n-1}a_{l}^{n}a_{\gamma}^{n}a_{\nu}^{\mu}c(dx_l)c(dx_\mu)c(dx_n)c(dx_\nu)c(dx_{\gamma})\\
&+\frac{i}{16(\xi_n-i)^2}h'(0)\sum_{l,\gamma,\mu=1}^{n}\sum_{q,\nu=1}^{n-1}\xi_{q}a_{l}^{q}a_{\gamma}^{n}a_{\nu}^{\mu}c(dx_l)c(dx_\mu)c(dx_n)c(dx_\nu)c(dx_{\gamma})\nonumber\\
&+\frac{i}{16(\xi_n-i)^2}h'(0)\sum_{l,\gamma,\mu=1}^{n}\sum_{\alpha,\nu=1}^{n-1}\xi_{\alpha}a_{l}^{n}a_{\gamma}^{\alpha}a_{\nu}^{\mu}c(dx_l)c(dx_\mu)c(dx_n)c(dx_\nu)c(dx_{\gamma})\nonumber\\
&+\frac{i\xi_n+2}{16(\xi_n-i)^2}h'(0)\sum_{l,\gamma,\mu=1}^{n}\sum_{q,\alpha,\nu=1}^{n-1}\xi_{q}\xi_{\alpha}a_{l}^{q}a_{\gamma}^{\alpha}a_{\nu}^{\mu}c(dx_l)c(dx_\mu)c(dx_n)c(dx_\nu)c(dx_{\gamma}).\nonumber
\end{align}
Accordingly, we have
\begin{align}
&{\rm trace} [\pi^+_{\xi_n}(A_1(x_0)) \times \partial_{\xi_n}\sigma_{-3}({{D}_{J}}^{-3})](x_0)|_{|\xi'|=1}\\
&=\frac{\xi _n^2}{4 \left(\xi _n-i\right)^5 \left(\xi _n+i\right)^3}h'(0)\sum_{l,\gamma,\mu,\beta=1}^{n}\sum_{\nu,i=1}^{n-1}{\rm tr}[\xi_{i}a_{l}^{n}a_{\gamma}^{n}a_{\nu}^{\mu}a_{\beta}^{i}c(dx_{l})c(dx_{\mu})c(dx_{n})c(dx_{\nu})c(dx_{\gamma})c(dx_{\beta})]\nonumber\\
&+\frac{\xi _n \left(3 \xi _n^2-1\right)}{16 \left(\xi _n-i\right)^5 \left(\xi _n+i\right)^3}h'(0)\sum_{l,\gamma,\mu,\beta=1}^{n}\sum_{\nu=1}^{n-1}{\rm tr}[a_{l}^{n}a_{\gamma}^{n}a_{\nu}^{\mu}a_{\beta}^{n}c(dx_{l})c(dx_{\mu})c(dx_{n})c(dx_{\nu})c(dx_{\gamma})c(dx_{\beta})]\nonumber\\
&+\frac{\xi _n}{4 \left(\xi _n-i\right)^5 \left(\xi _n+i\right)^3}h'(0)\sum_{l,\gamma,\mu,\beta=1}^{n}\sum_{q,\nu,i=1}^{n-1}{\rm tr}[\xi_{q}\xi_{i}a_{l}^{q}a_{\gamma}^{n}a_{\nu}^{\mu}a_{\beta}^{i}c(dx_{l})c(dx_{\mu})c(dx_{n})c(dx_{\nu})c(dx_{\gamma})c(dx_{\beta})]\nonumber\\
&+\frac{3 \xi _n^2-1}{16 \left(\xi _n-i\right)^5 \left(\xi _n+i\right)^3}h'(0)\sum_{l,\gamma,\mu,\beta=1}^{n}\sum_{q,\nu=1}^{n-1}{\rm tr}[\xi_{q}a_{l}^{q}a_{\gamma}^{n}a_{\nu}^{\mu}a_{\beta}^{n}c(dx_{l})c(dx_{\mu})c(dx_{n})c(dx_{\nu})c(dx_{\gamma})c(dx_{\beta})]\nonumber\\
&+\frac{\xi _n}{4 \left(\xi _n-i\right)^5 \left(\xi _n+i\right)^3}h'(0)\sum_{l,\gamma,\mu,\beta=1}^{n}\sum_{\alpha,\nu,i=1}^{n-1}{\rm tr}[\xi_{\alpha}\xi_{i}a_{l}^{n}a_{\gamma}^{\alpha}a_{\nu}^{\mu}a_{\beta}^{i}c(dx_{l})c(dx_{\mu})c(dx_{n})c(dx_{\nu})c(dx_{\gamma})c(dx_{\beta})]\nonumber\\
&+\frac{3 \xi _n^2-1}{16 \left(\xi _n-i\right)^5 \left(\xi _n+i\right)^3}h'(0)\sum_{l,\gamma,\mu,\beta=1}^{n}\sum_{\alpha,\nu=1}^{n-1}{\rm tr}[\xi_{\alpha}a_{l}^{n}a_{\gamma}^{\alpha}a_{\nu}^{\mu}a_{\beta}^{n}c(dx_{l})c(dx_{\mu})c(dx_{n})c(dx_{\nu})c(dx_{\gamma})c(dx_{\beta})]\nonumber\\
&+\frac{\xi _n \left(\xi _n-2 i\right)}{4 \left(\xi _n-i\right)^5 \left(\xi _n+i\right)^3}h'(0)\sum_{l,\gamma,\mu,\beta=1}^{n}\sum_{q,\alpha,\nu,i=1}^{n-1}{\rm tr}[\xi_{q}\xi_{\alpha}\xi_{i}a_{l}^{q}a_{\gamma}^{\alpha}a_{\nu}^{\mu}a_{\beta}^{i}c(dx_{l})c(dx_{\mu})c(dx_{n})c(dx_{\nu})c(dx_{\gamma})c(dx_{\beta})]\nonumber
\end{align}
\begin{align}
&+\frac{\left(\xi _n-2 i\right) \left(3 \xi _n^2-1\right)}{16 \left(\xi _n-i\right)^5 \left(\xi _n+i\right)^3}h'(0)\sum_{l,\gamma,\mu,\beta=1}^{n}\sum_{q,\alpha,\nu=1}^{n-1}{\rm tr}[\xi_{q}\xi_{\alpha}a_{l}^{q}a_{\gamma}^{\alpha}a_{\nu}^{\mu}a_{\beta}^{n}c(dx_{l})c(dx_{\mu})c(dx_{n})c(dx_{\nu})c(dx_{\gamma})c(dx_{\beta})].\nonumber
\end{align}
By $\int_{|\xi'|=1}{\{\xi_{i_1}\cdot\cdot\cdot\xi_{i_{2d+1}}}\}\sigma(\xi')=0$ and $\int_{|\xi'|=1}\xi_i\xi_j=\frac{\pi^{3}}{6}\delta_i^j,$ then we have
\begin{align}
&-i\int_{|\xi'|=1}\int^{+\infty}_{-\infty}{\rm trace} [\pi^+_{\xi_n}(A_1(x_0)) \times \partial_{\xi_n}\sigma_{-3}({{D}_{J}}^{-3})](x_0)d\xi_n\sigma(\xi')dx'\\
&=-i\Omega_4\int_{\Gamma^{+}}\frac{\xi _n \left(3 \xi _n^2-1\right)}{16 \left(\xi _n-i\right)^5 \left(\xi _n+i\right)^3}h'(0)\sum_{l=1}^{n}\sum_{\nu=1}^{n-1}(-(a_{\nu}^{n})^2(a_{l}^{n})^2+(a_{l}^{n})^2a_{\nu}^{\nu}a_{n}^{n}){\rm tr}[\texttt{id}]d\xi_ndx'\nonumber\\
&-i(\frac{16\pi^{2}}{15})\Omega_4\int_{\Gamma^{+}}\frac{\xi _n}{4 \left(\xi _n-i\right)^5 \left(\xi _n+i\right)^3}h'(0)\sum_{l=1}^{n}\sum_{\nu,i=1}^{n-1}(-(a_{\nu}^{n})^2(a_{l}^{i})^2+(a_{l}^{i})^2a_{\nu}^{\nu}a_{n}^{n}){\rm tr}[\texttt{id}]d\xi_ndx'\nonumber\\
&-i(\frac{8\pi^{2}}{15})\Omega_4\int_{\Gamma^{+}}\frac{\left(\xi _n-2 i\right) \left(3 \xi _n^2-1\right)}{16 \left(\xi _n-i\right)^5 \left(\xi _n+i\right)^3}h'(0)\sum_{l=1}^{n}\sum_{\nu,i=1}^{n-1}((a_{\nu}^{n})^2(a_{l}^{i})^2-(a_{l}^{i})^2a_{\nu}^{\nu}a_{n}^{n}-2a_{\nu}^{i}a_{l}^{i}a_{\nu}^{n}a_{l}^{n}+2a_{i}^{i}a_{\nu}^{\nu})\nonumber\\
&{\rm tr}[\texttt{id}]d\xi_ndx'.\nonumber
\end{align}
Hence
\begin{align}
&-i\int_{|\xi'|=1}\int^{+\infty}_{-\infty}{\rm trace} [\pi^+_{\xi_n}(A_1(x_0)) \times \partial_{\xi_n}\sigma_{-3}({{D}_{J}}^{-3})](x_0)d\xi_n\sigma(\xi')dx'\\
&=\sum_{l=1}^{n}\sum_{\nu=1}^{n-1}(-(a_{\nu}^{n})^2(a_{l}^{n})^2+(a_{l}^{n})^2a_{\nu}^{\nu}a_{n}^{n}){\rm tr}[\texttt{id}]\Omega_4h'(0)(-i)\frac{2\pi i}{4!}\Big[\frac{\xi _n \left(3 \xi _n^2-1\right)}{16 \left(\xi _n+i\right)^3}\Big]^{(4)}\bigg|_{\xi_n=i}dx'\nonumber\\
&+\sum_{l=1}^{n}\sum_{\nu,i=1}^{n-1}(-(a_{\nu}^{n})^2(a_{l}^{i})^2+(a_{l}^{i})^2a_{\nu}^{\nu}a_{n}^{n}){\rm tr}[\texttt{id}]\Omega_4h'(0)(-i)(\frac{16\pi^{2}}{15})\frac{2\pi i}{4!}\Big[\frac{\xi _n}{4 \left(\xi _n+i\right)^3}\Big]^{(4)}\bigg|_{\xi_n=i}dx'\nonumber\\
&+\sum_{l=1}^{n}\sum_{\nu,i=1}^{n-1}((a_{\nu}^{n})^2(a_{l}^{i})^2-(a_{l}^{i})^2a_{\nu}^{\nu}a_{n}^{n}-2a_{\nu}^{i}a_{l}^{i}a_{\nu}^{n}a_{l}^{n}+2a_{i}^{i}a_{\nu}^{\nu}){\rm tr}[\texttt{id}]\Omega_4h'(0)(-i)(\frac{8\pi^{2}}{15})\frac{2\pi i}{4!}\nonumber\\
&\Big[\frac{\left(\xi _n-2 i\right) \left(3 \xi _n^2-1\right)}{16 \left(\xi _n+i\right)^3}\Big]^{(4)}\bigg|_{\xi_n=i}dx'\nonumber\end{align}
\begin{align}
&=\sum_{l=1}^{n}\sum_{\nu=1}^{n-1}(-(a_{\nu}^{n})^2(a_{l}^{n})^2+(a_{l}^{n})^2a_{\nu}^{\nu}a_{n}^{n}){\rm tr}[\texttt{id}]\Omega_4h'(0)(\frac{\pi }{256})dx'\nonumber\\
&+\sum_{l=1}^{n}\sum_{\nu,i=1}^{n-1}(-(a_{\nu}^{n})^2(a_{l}^{i})^2+(a_{l}^{i})^2a_{\nu}^{\nu}a_{n}^{n}){\rm tr}[\texttt{id}]\Omega_4h'(0)(\frac{\pi ^3}{48})dx'\nonumber\\
&+\sum_{l=1}^{n}\sum_{\nu,i=1}^{n-1}((a_{\nu}^{n})^2(a_{l}^{i})^2-(a_{l}^{i})^2a_{\nu}^{\nu}a_{n}^{n}-2a_{\nu}^{i}a_{l}^{i}a_{\nu}^{n}a_{l}^{n}+2a_{i}^{i}a_{\nu}^{\nu}){\rm tr}[\texttt{id}]\Omega_4h'(0)(-\frac{\pi ^3}{96})dx'.\nonumber
\end{align}
\cite{LW2} also shown that
\begin{align}
\pi^+_{\xi_n}(A_2(x_0))
&=-\frac{i\xi_n}{4(\xi_n-i)^2}\sum_{l,j,h,y=1}^{n}a_{l}^{n}a_{y}^{j}\partial_{x_j}(a_{h}^{n})c(dx_{l})c(dx_{y})c(dx_{h})\nonumber\\
&-\frac{i}{4(\xi_n-i)^2}\sum_{l,j,h,y=1}^{n}\sum_{q=1}^{n-1}\xi_{q}a_{l}^{q}a_{y}^{j}\partial_{x_j}(a_{h}^{n})c(dx_{l})c(dx_{y})c(dx_{h})\nonumber\\
&-\frac{i}{4(\xi_n-i)^2}\sum_{l,j,h,y=1}^{n}\sum_{p=1}^{n-1}\xi_{p}a_{l}^{n}a_{y}^{j}\partial_{x_j}(a_{h}^{p})c(dx_{l})c(dx_{y})c(dx_{h})\nonumber\\
&-\frac{i\xi_n+2}{4(\xi_n-i)^2}\sum_{l,j,h,y=1}^{n}\sum_{q,p=1}^{n-1}\xi_{q}\xi_{p}a_{l}^{q}a_{y}^{j}\partial_{x_j}(a_{h}^{p})c(dx_{l})c(dx_{y})c(dx_{h})\nonumber\\
&-\frac{i\xi_n}{8(\xi_n-i)^2}h'(0)\sum_{l,z=1}^{n}\sum_{h=1}^{n-1}a_{l}^{n}a_{h}^{n}a_{z}^{n}c(dx_{l})c(dx_{z})c(dx_{h})\nonumber\\
&-\frac{i}{8(\xi_n-i)^2}h'(0)\sum_{l,z=1}^{n}\sum_{q,h=1}^{n-1}\xi_{q}a_{l}^{q}a_{h}^{n}a_{z}^{n}c(dx_{l})c(dx_{z})c(dx_{h})\nonumber\\
&-\frac{i}{8(\xi_n-i)^2}h'(0)\sum_{l,z=1}^{n}\sum_{p,h=1}^{n-1}\xi_{p}a_{l}^{n}a_{h}^{p}a_{z}^{n}c(dx_{l})c(dx_{z})c(dx_{h})\nonumber\\
&-\frac{i\xi_n+2}{8(\xi_n-i)^2}h'(0)\sum_{l,z=1}^{n}\sum_{q,p,h=1}^{n-1}\xi_{q}\xi_{p}a_{l}^{q}a_{h}^{p}a_{z}^{n}c(dx_{l})c(dx_{z})c(dx_{h}).\nonumber
\end{align}
It follows immediately that
\begin{align}
&{\rm trace} [\pi^+_{\xi_n}(A_2(x_0)) \times \partial_{\xi_n}\sigma_{-3}({{D}_{J}}^{-3})](x_0)|_{|\xi'|=1}\\
&=-\frac{\xi _n^2}{\left(\xi _n-i\right)^5 \left(\xi _n+i\right)^3}\sum_{l,j,h,y,\beta=1}^{n}\sum_{i=1}^{n-1}{\rm tr}[\xi_{i}a_{l}^{n}a_{y}^{j}a_{\beta}^{i}\partial_{x_j}(a_{h}^{n})c(dx_{l})c(dx_{y})c(dx_{h})c(dx_{\beta})]\nonumber\end{align}
\begin{align}
&+\frac{\xi _n-3 \xi _n^3}{4 \left(\xi _n-i\right)^5 \left(\xi _n+i\right)^3}\sum_{l,j,h,y,\beta=1}^{n}{\rm tr}[a_{l}^{n}a_{y}^{j}a_{\beta}^{n}\partial_{x_j}(a_{h}^{n})c(dx_{l})c(dx_{y})c(dx_{h})c(dx_{\beta})]\nonumber\\
&-\frac{\xi _n}{\left(\xi _n-i\right)^5 \left(\xi _n+i\right)^3}\sum_{l,j,h,y,\beta=1}^{n}\sum_{q,i=1}^{n-1}{\rm tr}[\xi_{q}\xi_{i}a_{l}^{q}a_{y}^{j}a_{\beta}^{i}\partial_{x_j}(a_{h}^{n})c(dx_{l})c(dx_{y})c(dx_{h})c(dx_{\beta})]\nonumber\\
&+\frac{1-3 \xi _n^2}{4 \left(\xi _n-i\right)^5 \left(\xi _n+i\right)^3}\sum_{l,j,h,y,\beta=1}^{n}\sum_{q=1}^{n-1}{\rm tr}[\xi_{q}a_{l}^{q}a_{y}^{j}a_{\beta}^{n}\partial_{x_j}(a_{h}^{n})c(dx_{l})c(dx_{y})c(dx_{h})c(dx_{\beta})]\nonumber\\
&-\frac{\xi _n}{\left(\xi _n-i\right)^5 \left(\xi _n+i\right)^3}\sum_{l,j,h,y,\beta=1}^{n}\sum_{p,i=1}^{n-1}{\rm tr}[\xi_{p}\xi_{i}a_{l}^{n}a_{y}^{j}a_{\beta}^{i}\partial_{x_j}(a_{h}^{p})c(dx_{l})c(dx_{y})c(dx_{h})c(dx_{\beta})]\nonumber\\
&+\frac{1-3 \xi _n^2}{4 \left(\xi _n-i\right)^5 \left(\xi _n+i\right)^3}\sum_{l,j,h,y,\beta=1}^{n}\sum_{p=1}^{n-1}{\rm tr}[\xi_{p}a_{l}^{n}a_{y}^{j}a_{\beta}^{n}\partial_{x_j}(a_{h}^{p})c(dx_{l})c(dx_{y})c(dx_{h})c(dx_{\beta})]\nonumber\\
&-\frac{\xi _n \left(\xi _n-2 i\right)}{\left(\xi _n-i\right)^5 \left(\xi _n+i\right)^3}\sum_{l,j,h,y,\beta=1}^{n}\sum_{q,p,i=1}^{n-1}{\rm tr}[\xi_{q}\xi_{p}\xi_{i}a_{l}^{q}a_{y}^{j}a_{\beta}^{i}\partial_{x_j}(a_{h}^{p})c(dx_{l})c(dx_{y})c(dx_{h})c(dx_{\beta})]\nonumber\\
&-\frac{\left(\xi _n-2 i\right) \left(3 \xi _n^2-1\right)}{4 \left(\xi _n-i\right)^5 \left(\xi _n+i\right)^3}\sum_{l,j,h,y,\beta=1}^{n}\sum_{q,p=1}^{n-1}{\rm tr}[\xi_{q}\xi_{p}a_{l}^{q}a_{y}^{j}a_{\beta}^{n}\partial_{x_j}(a_{h}^{p})c(dx_{l})c(dx_{y})c(dx_{h})c(dx_{\beta})]\nonumber\\
&-\frac{\xi _n^2}{2 \left(\xi _n-i\right)^5 \left(\xi _n+i\right)^3}h'(0)\sum_{l,z,\beta=1}^{n}\sum_{h,i=1}^{n-1}{\rm tr}[\xi_{i}a_{l}^{n}a_{h}^{n}a_{z}^{n}a_{\beta}^{i}c(dx_{l})c(dx_{z})c(dx_{h})c(dx_{\beta})]\nonumber\\
&+\frac{\xi _n-3 \xi _n^3}{8 \left(\xi _n-i\right)^5 \left(\xi _n+i\right)^3}h'(0)\sum_{l,z,\beta=1}^{n}\sum_{h=1}^{n-1}{\rm tr}[a_{l}^{n}a_{h}^{n}a_{z}^{n}a_{\beta}^{n}c(dx_{l})c(dx_{z})c(dx_{h})c(dx_{\beta})]\nonumber\\
&-\frac{\xi _n}{2 \left(\xi _n-i\right)^5 \left(\xi _n+i\right)^3}h'(0)\sum_{l,z,\beta=1}^{n}\sum_{q,h,i=1}^{n-1}{\rm tr}[\xi_{q}\xi_{i}a_{l}^{q}a_{h}^{n}a_{z}^{n}a_{\beta}^{i}c(dx_{l})c(dx_{z})c(dx_{h})c(dx_{\beta})]\nonumber\\
&+\frac{1-3 \xi _n^2}{8 \left(\xi _n-i\right)^5 \left(\xi _n+i\right)^3}h'(0)\sum_{l,z,\beta=1}^{n}\sum_{q,h=1}^{n-1}{\rm tr}[\xi_{q}a_{l}^{q}a_{h}^{n}a_{z}^{n}a_{\beta}^{n}c(dx_{l})c(dx_{z})c(dx_{h})c(dx_{\beta})]\nonumber\\
&-\frac{\xi _n}{2 \left(\xi _n-i\right)^5 \left(\xi _n+i\right)^3}h'(0)\sum_{l,z,\beta=1}^{n}\sum_{p,h,i=1}^{n-1}{\rm tr}[\xi_{p}\xi_{i}a_{l}^{n}a_{h}^{p}a_{z}^{n}a_{\beta}^{i}c(dx_{l})c(dx_{z})c(dx_{h})c(dx_{\beta})]\nonumber\\
&+\frac{1-3 \xi _n^2}{8 \left(\xi _n-i\right)^5 \left(\xi _n+i\right)^3}h'(0)\sum_{l,z,\beta=1}^{n}\sum_{p,h=1}^{n-1}{\rm tr}[\xi_{p}a_{l}^{n}a_{h}^{p}a_{z}^{n}a_{\beta}^{n}c(dx_{l})c(dx_{z})c(dx_{h})c(dx_{\beta})]\nonumber\\
&-\frac{\xi _n \left(\xi _n-2 i\right)}{2 \left(\xi _n-i\right)^5 \left(\xi _n+i\right)^3}h'(0)\sum_{l,z,\beta=1}^{n}\sum_{q,p,h,i=1}^{n-1}{\rm tr}[\xi_{q}\xi_{p}\xi_{i}a_{l}^{q}a_{h}^{p}a_{z}^{n}a_{\beta}^{i}c(dx_{l})c(dx_{z})c(dx_{h})c(dx_{\beta})]\nonumber\\
&-\frac{\left(\xi _n-2 i\right) \left(3 \xi _n^2-1\right)}{8 \left(\xi _n-i\right)^5 \left(\xi _n+i\right)^3}h'(0)\sum_{l,z,\beta=1}^{n}\sum_{q,p,h=1}^{n-1}{\rm tr}[\xi_{q}\xi_{p}a_{l}^{q}a_{h}^{p}a_{z}^{n}a_{\beta}^{n}c(dx_{l})c(dx_{z})c(dx_{h})c(dx_{\beta})].\nonumber
\end{align}
Similarly, we have
\begin{align}
&-i\int_{|\xi'|=1}\int^{+\infty}_{-\infty}{\rm trace} [\pi^+_{\xi_n}(A_2(x_0)) \times \partial_{\xi_n}\sigma_{-3}({{D}_{J}}^{-3})](x_0)d\xi_n\sigma(\xi')dx'\\
&=-i\Omega_4\int_{\Gamma^{+}}\frac{\xi _n-3 \xi _n^3}{4 \left(\xi _n-i\right)^5 \left(\xi _n+i\right)^3}\sum_{l,j,\beta=1}^{n}\left((a_{\beta}^{n})^2a_{l}^{j}\partial_{x_j}(a_{l}^{n})-a_{l}^{n}a_{\beta}^{j}a_{\beta}^{n}\partial_{x_j}(a_{l}^{n})+a_{l}^{n}a_{l}^{j}a_{\beta}^{n}\partial_{x_j}(a_{\beta}^{n})\right)\nonumber\\
&{\rm tr}[\texttt{id}]d\xi_ndx'\nonumber\\
&-i(\frac{8\pi^{2}}{15})\Omega_4\int_{\Gamma^{+}}-\frac{\xi _n}{\left(\xi _n-i\right)^5 \left(\xi _n+i\right)^3}\sum_{l,j,\beta=1}^{n}\sum_{i=1}^{n-1}\left((a_{\beta}^{i})^2a_{l}^{j}\partial_{x_j}(a_{l}^{n})-a_{l}^{i}a_{\beta}^{j}a_{\beta}^{i}\partial_{x_j}(a_{l}^{n})+a_{l}^{i}a_{l}^{j}a_{\beta}^{i}\partial_{x_j}(a_{\beta}^{n})\right)\nonumber\\
&{\rm tr}[\texttt{id}]d\xi_ndx'\nonumber\\
&-i(\frac{8\pi^{2}}{15})\Omega_4\int_{\Gamma^{+}}-\frac{\xi _n}{\left(\xi _n-i\right)^5 \left(\xi _n+i\right)^3}\sum_{l,j,\beta=1}^{n}\sum_{i=1}^{n-1}\left(a_{\beta}^{n}a_{l}^{j}a_{\beta}^{i}\partial_{x_j}(a_{l}^{i})-a_{l}^{n}a_{\beta}^{j}a_{\beta}^{i}\partial_{x_j}(a_{l}^{i})+a_{l}^{n}a_{l}^{j}a_{\beta}^{i}\partial_{x_j}(a_{\beta}^{i})\right)\nonumber\\
&{\rm tr}[\texttt{id}]d\xi_ndx'\nonumber\\
&-i(\frac{8\pi^{2}}{15})\Omega_4\int_{\Gamma^{+}}-\frac{\left(\xi _n-2 i\right) \left(3 \xi _n^2-1\right)}{4 \left(\xi _n-i\right)^5 \left(\xi _n+i\right)^3}
\sum_{l,j,\beta=1}^{n}\sum_{i=1}^{n-1}\left(a_{\beta}^{i}a_{l}^{j}a_{\beta}^{n}\partial_{x_j}(a_{l}^{i})-a_{l}^{i}a_{\beta}^{j}a_{\beta}^{n}\partial_{x_j}(a_{l}^{i})+a_{l}^{i}a_{l}^{j}a_{\beta}^{n}\partial_{x_j}(a_{\beta}^{i})\right)\nonumber\\
&{\rm tr}[\texttt{id}]d\xi_ndx'\nonumber\\
&-i\Omega_4\int_{\Gamma^{+}}\frac{\xi _n-3 \xi _n^3}{8 \left(\xi _n-i\right)^5 \left(\xi _n+i\right)^3}
h'(0)\sum_{l=1}^{n}\sum_{\nu=1}^{n-1}(a_{\nu}^{n})^2(a_{l}^{n})^2{\rm tr}[\texttt{id}]d\xi_ndx'\nonumber\\
&-i(\frac{8\pi^{2}}{15})\Omega_4\int_{\Gamma^{+}}-\frac{\xi _n}{2 \left(\xi _n-i\right)^5 \left(\xi _n+i\right)^3}
h'(0)\sum_{l=1}^{n}\sum_{\nu,i=1}^{n-1}(a_{\nu}^{n})^2(a_{l}^{i})^2{\rm tr}[\texttt{id}]d\xi_ndx'\nonumber\\
&-i(\frac{8\pi^{2}}{15})\Omega_4\int_{\Gamma^{+}}-\frac{\xi _n}{2 \left(\xi _n-i\right)^5 \left(\xi _n+i\right)^3}
h'(0)\sum_{l=1}^{n}\sum_{\nu,i=1}^{n-1}(a_{\nu}^{i})^2(a_{l}^{n})^2{\rm tr}[\texttt{id}]d\xi_ndx'\nonumber\\
&-i(\frac{8\pi^{2}}{15})\Omega_4\int_{\Gamma^{+}}-\frac{\left(\xi _n-2 i\right) \left(3 \xi _n^2-1\right)}{8 \left(\xi _n-i\right)^5 \left(\xi _n+i\right)^3}
h'(0)\sum_{l=1}^{n}\sum_{\nu,i=1}^{n-1}\left(2a_{\nu}^{i}a_{l}^{i}a_{\nu}^{n}a_{l}^{n}-(a_{\nu}^{i})^2(a_{l}^{n})^2\right){\rm tr}[\texttt{id}]d\xi_ndx'.\nonumber
\end{align}
A simple calculation shows that
\begin{align}
&-i\int_{|\xi'|=1}\int^{+\infty}_{-\infty}{\rm trace} [\pi^+_{\xi_n}(A_2(x_0)) \times \partial_{\xi_n}\sigma_{-3}({{D}_{J}}^{-3})](x_0)d\xi_n\sigma(\xi')dx'\\
&=\sum_{l,j,\beta=1}^{n}\left((a_{\beta}^{n})^2a_{l}^{j}\partial_{x_j}(a_{l}^{n})-a_{l}^{n}a_{\beta}^{j}a_{\beta}^{n}\partial_{x_j}(a_{l}^{n})+a_{l}^{n}a_{l}^{j}a_{\beta}^{n}\partial_{x_j}(a_{\beta}^{n})\right){\rm tr}[\texttt{id}]\Omega_4(-\frac{\pi }{64})dx'\nonumber
\end{align}
\begin{align}
&+\sum_{l,j,\beta=1}^{n}\sum_{i=1}^{n-1}\left((a_{\beta}^{i})^2a_{l}^{j}\partial_{x_j}(a_{l}^{n})-a_{l}^{i}a_{\beta}^{j}a_{\beta}^{i}\partial_{x_j}(a_{l}^{n})+a_{l}^{i}a_{l}^{j}a_{\beta}^{i}\partial_{x_j}(a_{\beta}^{n})\right){\rm tr}[\texttt{id}]\Omega_4(-\frac{\pi ^3}{24} )dx'\nonumber\\
&+\sum_{l,j,\beta=1}^{n}\sum_{i=1}^{n-1}\left(a_{\beta}^{n}a_{l}^{j}a_{\beta}^{i}\partial_{x_j}(a_{l}^{i})-a_{l}^{n}a_{\beta}^{j}a_{\beta}^{i}\partial_{x_j}(a_{l}^{i})+a_{l}^{n}a_{l}^{j}a_{\beta}^{i}\partial_{x_j}(a_{\beta}^{i})\right){\rm tr}[\texttt{id}]\Omega_4(-\frac{\pi ^3}{24})dx'\nonumber\\
&+\sum_{l,j,\beta=1}^{n}\sum_{i=1}^{n-1}\left(a_{\beta}^{i}a_{l}^{j}a_{\beta}^{n}\partial_{x_j}(a_{l}^{i})-a_{l}^{i}a_{\beta}^{j}a_{\beta}^{n}\partial_{x_j}(a_{l}^{i})+a_{l}^{i}a_{l}^{j}a_{\beta}^{n}\partial_{x_j}(a_{\beta}^{i})\right){\rm tr}[\texttt{id}]\Omega_4(\frac{\pi ^3}{24})dx'\nonumber\\
&+\sum_{l=1}^{n}\sum_{\nu=1}^{n-1}(a_{\nu}^{n})^2(a_{l}^{n})^2{\rm tr}[\texttt{id}]\Omega_4h'(0)(-\frac{\pi }{128})dx'\nonumber\\
&+\sum_{l=1}^{n}\sum_{\nu,i=1}^{n-1}(a_{\nu}^{n})^2(a_{l}^{i})^2{\rm tr}[\texttt{id}]\Omega_4h'(0)(-\frac{\pi ^3}{48} )dx'\nonumber\\
&+\sum_{l=1}^{n}\sum_{\nu,i=1}^{n-1}(a_{\nu}^{i})^2(a_{l}^{n})^2{\rm tr}[\texttt{id}]\Omega_4h'(0)(-\frac{\pi ^3}{48} )dx'\nonumber\\
&+\sum_{l=1}^{n}\sum_{\nu,i=1}^{n-1}\left(2a_{\nu}^{i}a_{l}^{i}a_{\nu}^{n}a_{l}^{n}-(a_{\nu}^{i})^2(a_{l}^{n})^2\right){\rm tr}[\texttt{id}]\Omega_4h'(0)(\frac{\pi ^3}{48})dx'.\nonumber
\end{align}
Since
\begin{align}
-h'(0)\pi^+_{\xi_n}(A_3(x_0))
&=\frac{3\xi_n+i\xi_n^2}{16(\xi_n-i)^3}h'(0)\sum_{l,w,\gamma=1}^{n}a_{l}^{n}a_{w}^{n}a_{\gamma}^{n}c(dx_{l})c(dx_{w})c(dx_{\gamma})\\
&+\frac{i\xi_n+3}{16(\xi_n-i)^3}h'(0)\sum_{l,w,\gamma=1}^{n}\sum_{q=1}^{n-1}\xi_{q}a_{l}^{q}a_{w}^{n}a_{\gamma}^{n}c(dx_{l})c(dx_{w})c(dx_{\gamma})\nonumber\\
&+\frac{i\xi_n+3}{16(\xi_n-i)^3}h'(0)\sum_{l,w,\gamma=1}^{n}\sum_{\alpha=1}^{n-1}\xi_{\alpha}a_{l}^{n}a_{w}^{n}a_{\gamma}^{\alpha}c(dx_{l})c(dx_{w})c(dx_{\gamma})\nonumber\\
&+\frac{-8i+9\xi_n+3i\xi_n^2}{16(\xi_n-i)^3}h'(0)\sum_{l,w,\gamma=1}^{n}\sum_{q,\alpha=1}^{n-1}\xi_{q}\xi_{\alpha}a_{l}^{q}a_{w}^{n}a_{\gamma}^{\alpha}c(dx_{l})c(dx_{w})c(dx_{\gamma}),\nonumber
\end{align}
it is sufficient to show that
\begin{align}
&{\rm trace} [-h'(0)\pi^+_{\xi_n}(A_3(x_0)) \times \partial_{\xi_n}\sigma_{-3}({{D}_{J}}^{-3})](x_0)|_{|\xi'|=1}\\
&=\frac{\xi _n^2 \left(\xi _n-3 i\right)}{4 \left(\xi _n-i\right)^6 \left(\xi _n+i\right)^3}h'(0)\sum_{l,w,\gamma,\beta=1}^{n}\sum_{i=1}^{n-1}{\rm tr}[\xi_{i}a_{l}^{n}a_{w}^{n}a_{\gamma}^{n}a_{\beta}^{i}c(dx_{l})c(dx_{w})c(dx_{\gamma})c(dx_{\beta})]\nonumber\end{align}
\begin{align}
&+\frac{\xi _n \left(\xi _n-3 i\right) \left(3 \xi _n^2-1\right)}{16 \left(\xi _n-i\right)^6 \left(\xi _n+i\right)^3}h'(0)\sum_{l,w,\gamma,\beta=1}^{n}{\rm tr}[a_{l}^{n}a_{w}^{n}a_{\gamma}^{n}a_{\beta}^{n}c(dx_{l})c(dx_{w})c(dx_{\gamma})c(dx_{\beta})]\nonumber\\
&+\frac{\xi _n \left(\xi _n-3 i\right)}{4 \left(\xi _n-i\right)^6 \left(\xi _n+i\right)^3}h'(0)\sum_{l,w,\gamma,\beta=1}^{n}\sum_{q,i=1}^{n-1}{\rm tr}[\xi_{q}\xi_{i}a_{l}^{q}a_{w}^{n}a_{\gamma}^{n}a_{\beta}^{i}c(dx_{l})c(dx_{w})c(dx_{\gamma})c(dx_{\beta})]\nonumber\\
&+\frac{\left(\xi _n-3 i\right) \left(3 \xi _n^2-1\right)}{16 \left(\xi _n-i\right)^6 \left(\xi _n+i\right)^3}h'(0)\sum_{l,w,\gamma,\beta=1}^{n}\sum_{q=1}^{n-1}{\rm tr}[\xi_{q}a_{l}^{q}a_{w}^{n}a_{\gamma}^{n}a_{\beta}^{n}c(dx_{l})c(dx_{w})c(dx_{\gamma})c(dx_{\beta})]\nonumber\\
&+\frac{\xi _n \left(\xi _n-3 i\right)}{4 \left(\xi _n-i\right)^6 \left(\xi _n+i\right)^3}h'(0)\sum_{l,w,\gamma,\beta=1}^{n}\sum_{\alpha,i=1}^{n-1}{\rm tr}[\xi_{\alpha}\xi_{i}a_{l}^{n}a_{w}^{n}a_{\gamma}^{\alpha}a_{\beta}^{i}c(dx_{l})c(dx_{w})c(dx_{\gamma})c(dx_{\beta})]\nonumber\\
&+\frac{\left(\xi _n-3 i\right) \left(3 \xi _n^2-1\right)}{16 \left(\xi _n-i\right)^6 \left(\xi _n+i\right)^3}h'(0)\sum_{l,w,\gamma,\beta=1}^{n}\sum_{\alpha=1}^{n-1}{\rm tr}[\xi_{\alpha}a_{l}^{n}a_{w}^{n}a_{\gamma}^{\alpha}a_{\beta}^{n}c(dx_{l})c(dx_{w})c(dx_{\gamma})c(dx_{\beta})]\nonumber\\
&+\frac{\xi _n \left(3 \xi _n^2-9 i \xi _n-8\right)}{4 \left(\xi _n-i\right){}^6 \left(\xi _n+i\right){}^3}h'(0)\sum_{l,w,\gamma,\beta=1}^{n}\sum_{q,\alpha,i=1}^{n-1}{\rm tr}[\xi_{q}\xi_{\alpha}\xi_{i}a_{l}^{q}a_{w}^{n}a_{\gamma}^{\alpha}a_{\beta}^{i}c(dx_{l})c(dx_{w})c(dx_{\gamma})c(dx_{\beta})]\nonumber\\
&+\frac{\left(3 \xi _n^2-1\right) \left(3 \xi _n^2-9 i \xi _n-8\right)}{16 \left(\xi _n-i\right){}^6 \left(\xi _n+i\right){}^3}h'(0)\sum_{l,w,\gamma,\beta=1}^{n}\sum_{q,\alpha=1}^{n-1}{\rm tr}[\xi_{q}\xi_{\alpha}a_{l}^{q}a_{w}^{n}a_{\gamma}^{\alpha}a_{\beta}^{n}c(dx_{l})c(dx_{w})c(dx_{\gamma})c(dx_{\beta})],\nonumber
\end{align}
then, we have
\begin{align}
&-i\int_{|\xi'|=1}\int^{+\infty}_{-\infty}{\rm trace} [-h'(0)\pi^+_{\xi_n}(A_3(x_0)) \times \partial_{\xi_n}\sigma_{-3}({{D}_{J}}^{-3})](x_0)d\xi_n\sigma(\xi')dx'\\
&=-i\Omega_4\int_{\Gamma^{+}}\frac{\xi _n \left(\xi _n-3 i\right) \left(3 \xi _n^2-1\right)}{16 \left(\xi _n-i\right)^6 \left(\xi _n+i\right)^3}h'(0)\sum_{\beta,l=1}^{n}(a_{\beta}^{n})^2(a_{l}^{n})^2{\rm tr}[\texttt{id}]d\xi_ndx'\nonumber\\
&-i(\frac{16 \pi^{2}}{15})\Omega_4\int_{\Gamma^{+}}\frac{\xi _n \left(\xi _n-3 i\right)}{4 \left(\xi _n-i\right)^6 \left(\xi _n+i\right)^3}h'(0)\sum_{\beta,l=1}^{n}\sum_{i=1}^{n-1}(a_{\beta}^{i})^2(a_{l}^{n})^2{\rm tr}[\texttt{id}]d\xi_ndx'\nonumber\\
&-i(\frac{8 \pi^{2}}{15})\Omega_4\int_{\Gamma^{+}}\frac{\left(3 \xi _n^2-1\right) \left(3 \xi _n^2-9 i \xi _n-8\right)}{16 \left(\xi _n-i\right){}^6 \left(\xi _n+i\right){}^3}h'(0)\sum_{\beta,l=1}^{n}\sum_{i=1}^{n-1}\left(2a_{l}^{i}a_{\beta}^{i}a_{l}^{n}a_{\beta}^{n}-(a_{l}^{i})^2(a_{\beta}^{n})^2\right){\rm tr}[\texttt{id}]d\xi_ndx'\nonumber\\
&=\sum_{\beta,l=1}^{n}(a_{\beta}^{n})^2(a_{l}^{n})^2{\rm tr}[\texttt{id}]\Omega_4h'(0)(\frac{\pi }{256})dx'+\sum_{\beta,l=1}^{n}\sum_{i=1}^{n-1}(a_{\beta}^{i})^2(a_{l}^{n})^2{\rm tr}[\texttt{id}]\Omega_4h'(0)(\frac{7 \pi ^3}{120})dx'\nonumber\\
&+\sum_{\beta,l=1}^{n}\sum_{i=1}^{n-1}\left(2a_{l}^{i}a_{\beta}^{i}a_{l}^{n}a_{\beta}^{n}-(a_{l}^{i})^2(a_{\beta}^{n})^2\right){\rm tr}[\texttt{id}]\Omega_4h'(0)(-\frac{7 \pi ^3}{160} )dx'.\nonumber
\end{align}
On account of the above result,
\begin{align}
\Phi_5&=\sum_{l,\beta=1}^{n}\sum_{i=1}^{n-1}a_{l}^{i}a_{\beta}^{i}a_{l}^{n}a_{\beta}^{n}{\rm tr}[\texttt{id}]\Omega_4h'(0)(-\frac{7 \pi ^3}{80} )dx'
+\sum_{l=1}^{n}\sum_{\nu,i=1}^{n-1}a_{\nu}^{i}a_{l}^{i}a_{\nu}^{n}a_{l}^{n}{\rm tr}[\texttt{id}]\Omega_4h'(0)(\frac{\pi ^3}{16})dx'\nonumber\\
&+\sum_{l,j,\beta=1}^{n}(a_{\beta}^{n})^2a_{l}^{j}\partial_{x_j}(a_{l}^{n}){\rm tr}[\texttt{id}]\Omega_4(-\frac{\pi }{64})dx'
+\sum_{l,j,\beta=1}^{n}\sum_{i=1}^{n-1}(a_{\beta}^{i})^2a_{l}^{j}\partial_{x_j}(a_{l}^{n}){\rm tr}[\texttt{id}]\Omega_4(-\frac{\pi ^3}{24} )dx'\nonumber\\
&+\sum_{l=1}^{n}\sum_{i=1}^{n-1}(a_{l}^{n})^2a_{i}^{i}a_{n}^{n}{\rm tr}[\texttt{id}]\Omega_4h'(0)(\frac{\pi }{256})dx'
+\sum_{l=1}^{n}\sum_{\nu,i=1}^{n-1}(a_{l}^{i})^2a_{\nu}^{\nu}a_{n}^{n}{\rm tr}[\texttt{id}]\Omega_4h'(0)(\frac{\pi ^3}{32})dx'\nonumber\\
&+\sum_{l,j=1}^{n}a_{l}^{j}\partial_{x_j}(a_{l}^{n}){\rm tr}[\texttt{id}]\Omega_4(-\frac{\pi^{3} }{12})dx'
+\sum_{l,\beta=1}^{n}(a_{\beta}^{n})^2(a_{l}^{n})^2{\rm tr}[\texttt{id}]\Omega_4h'(0)(\frac{\pi }{256})dx'\nonumber\\
&+\sum_{l,\beta=1}^{n}\sum_{i=1}^{n-1}(a_{l}^{n})^2(a_{\beta}^{i})^2{\rm tr}[\texttt{id}]\Omega_4h'(0)(\frac{49 \pi ^3}{480})dx'
+\sum_{l=1}^{n}\sum_{i=1}^{n-1}(a_{i}^{n})^2(a_{l}^{n})^2{\rm tr}[\texttt{id}]\Omega_4h'(0)(-\frac{3 \pi }{256})dx'\nonumber\\
&+\sum_{l=1}^{n}\sum_{\nu,i=1}^{n-1}(a_{\nu}^{n})^2(a_{l}^{i})^2{\rm tr}[\texttt{id}]\Omega_4h'(0)(-\frac{5 \pi ^3}{96} )dx'
+\sum_{l=1}^{n}\sum_{\nu,i=1}^{n-1}(a_{l}^{n})^2(a_{\nu}^{i})^2{\rm tr}[\texttt{id}]\Omega_4h'(0)(-\frac{ \pi ^3}{24} )dx'\nonumber\\
&+\sum_{\nu,i=1}^{n-1}a_{i}^{i}a_{\nu}^{\nu}{\rm tr}[\texttt{id}]\Omega_4h'(0)(-\frac{\pi ^3}{48})dx'.\nonumber
\end{align}

In summary,
\begin{align}
\Phi
&=\Phi_1+\Phi_2+\Phi_3+\Phi_4+\Phi_5\\
&=\Big(\sum_{l=1}^{n}g^{M}(J(dx_{l}), (\nabla^{L}_{e_{l}}J)e_{n}){\rm tr}[\texttt{id}]
-\sum_{l=1}^{n}g^{M}(J(dx_{n}), (\nabla^{L}_{e_{l}}J)e_{l}){\rm tr}[\texttt{id}]\nonumber\\
&+\sum_{l=1}^{n}g^{M}(J(dx_{l}), (\nabla^{L}_{e_{n}}J)e_{l}){\rm tr}[\texttt{id}]\Big)
\Omega_4(-\frac{\pi ^3}{8})dx'
+\sum_{l,j=1}^{n}a_{l}^{j}\partial_{x_{j}}(a_{l}^{n}){\rm tr}[\texttt{id}]\Omega_4(\frac{1}{48} \pi  \left(10 \pi ^2-3\right))dx'\nonumber\\
&+\sum_{l,j,\beta=1}^{n}(a_{\beta}^{n})^2a_{l}^{j}\partial_{x_j}(a_{l}^{n}){\rm tr}[\texttt{id}]\Omega_4(-\frac{\pi }{64})dx'
+\sum_{l,j,\beta=1}^{n}\sum_{i=1}^{n-1}(a_{\beta}^{i})^2a_{l}^{j}\partial_{x_j}(a_{l}^{n}){\rm tr}[\texttt{id}]\Omega_4(-\frac{\pi ^3}{24} )dx'\nonumber\\
&+\sum_{l,\beta=1}^{n}\sum_{i=1}^{n-1}a_{l}^{i}a_{\beta}^{i}a_{l}^{n}a_{\beta}^{n}{\rm tr}[\texttt{id}]\Omega_4h'(0)(-\frac{7 \pi ^3}{80} )dx'
+\sum_{l=1}^{n}\sum_{\nu,i=1}^{n-1}a_{\nu}^{i}a_{l}^{i}a_{\nu}^{n}a_{l}^{n}{\rm tr}[\texttt{id}]\Omega_4h'(0)(\frac{\pi ^3}{16})dx'\nonumber\\
&+\sum_{l=1}^{n}\sum_{i=1}^{n-1}(a_{l}^{n})^2a_{i}^{i}a_{n}^{n}{\rm tr}[\texttt{id}]\Omega_4h'(0)(\frac{\pi }{256})dx'
+\sum_{l=1}^{n}\sum_{\nu,i=1}^{n-1}(a_{l}^{i})^2a_{\nu}^{\nu}a_{n}^{n}{\rm tr}[\texttt{id}]\Omega_4h'(0)(\frac{\pi ^3}{32})dx'\nonumber\end{align}
\begin{align}
&+\sum_{l,\beta=1}^{n}(a_{\beta}^{n})^2(a_{l}^{n})^2{\rm tr}[\texttt{id}]\Omega_4h'(0)(\frac{\pi }{256})dx'
+\sum_{l,\beta=1}^{n}\sum_{i=1}^{n-1}(a_{l}^{n})^2(a_{\beta}^{i})^2{\rm tr}[\texttt{id}]\Omega_4h'(0)(\frac{2 \pi ^3}{15})dx'\nonumber\\
&+\sum_{l=1}^{n}\sum_{i=1}^{n-1}(a_{i}^{n})^2(a_{l}^{n})^2{\rm tr}[\texttt{id}]\Omega_4h'(0)(-\frac{3 \pi}{256})dx'
+\sum_{l=1}^{n}\sum_{\nu,i=1}^{n-1}(a_{\nu}^{n})^2(a_{l}^{i})^2{\rm tr}[\texttt{id}]\Omega_4h'(0)(-\frac{5 \pi ^3}{96})dx'\nonumber\\
&+\sum_{l=1}^{n}\sum_{\nu,i=1}^{n-1}(a_{l}^{n})^2(a_{\nu}^{i})^2{\rm tr}[\texttt{id}]\Omega_4h'(0)(-\frac{\pi ^3}{24} )dx'
+\sum_{i=1}^{n-1}a_{n}^{n}a_{i}^{i}{\rm tr}[\texttt{id}]\Omega_4h'(0)(-\frac{3 \pi^3}{64})dx'\nonumber\\
&+\sum_{\nu,i=1}^{n-1}a_{i}^{i}a_{\nu}^{\nu}{\rm tr}[\texttt{id}]\Omega_4h'(0)(-\frac{\pi ^3}{48})dx'
+\sum_{l=1}^{n}(a_{l}^{n})^{2}{\rm tr}[\texttt{id}]\Omega_4h'(0)(\frac{1}{640} \pi  \left(384 \pi ^2+15\right))dx'\nonumber\\
&+\sum_{i=1}^{n-1}(a_{i}^{n})^{2}{\rm tr}[\texttt{id}]\Omega_4h'(0)(-\frac{5 \pi^3}{64})dx'
+\sum_{l=1}^{n}\sum_{i=1}^{n-1}(a_{l}^{i})^2{\rm tr}[\texttt{id}]\Omega_4h'(0)(\frac{7 \pi ^3}{240})dx'.\nonumber
\end{align}

(2) makes it obvious that
\begin{align}
\sum_{\beta=1}^{n}\sum_{i=1}^{n-1}a_{\beta}^{i}\partial_{x_i}(a_{\beta}^{n})=
\sum_{\beta=1}^{n}\langle\nabla_{J(e_{\beta})}^{L}(Je_{n}), e_{\beta}\rangle-\sum_{\beta=1}^{n}g^{M}\left(J(\frac{\partial}{\partial{x_{n}}}), \frac{\partial}{\partial{x_{n}}}\right)\langle\nabla_{J(e_{\beta})}^{L}(\frac{\partial}{\partial{x_{n}}}), e_{\beta}\rangle.
\end{align}

Combine (3.19) with (3.89), we obtain immediately the following theorem:
\begin{thm}
Let $M$ be a $6$-dimensional almost product Riemannian spin manifold with the boundary $\partial M$ and the metric
$g^M$ as above, ${{D}_{J}}$ be the $J$-twist of the Dirac operator on $\widetilde{M}$, then
\begin{align}
\widetilde{{\rm Wres}}[\pi^+{{D}_{J}}^{-1}\circ\pi^+{{D}_{J}}^{-3}]&=\int_{M}256\pi^{3}\Big(\sum_{i,j=1}^{6}R(J(e_{i}), J(e_{j}), e_{j}, e_{i})
-2\sum_{\nu,j=1}^{6}g^{M}(\nabla_{e_{j}}^{L}(J)e_{\nu}, (\nabla^{L}_{e_{\nu}}J)e_{j})\nonumber\\
&-2\sum_{\nu,j=1}^{6}g^{M}(J(e_{\nu}), (\nabla^{L}_{e_{j}}(\nabla^{L}_{e_{\nu}}(J)))e_{j}-(\nabla^{L}_{\nabla^{L}_{e_{j}}e_{\nu}}(J))e_{j})\nonumber\\
&-\sum_{\alpha,\nu,j=1}^{6}g^{M}(J(e_{\alpha}), (\nabla^{L}_{e_{\nu}}J)e_{j})g^{M}((\nabla^{L}_{e_{\alpha}}J)e_{j}, J(e_{\nu}))\nonumber\\
&-\sum_{\alpha,\nu,j=1}^{6}g^{M}(J(e_{\alpha}), (\nabla^{L}_{e_{\alpha}}J)e_{j})g^{M}(J(e_{\nu}), (\nabla^{L}_{e_{\nu}}J)e_{j})\nonumber\\
&+\sum_{\nu,j=1}^{6}g^{M}((\nabla^{L}_{e_{\nu}}J)e_{j}, (\nabla^{L}_{e_{\nu}}J)e_{j}))-\frac{1}{3}s\Big)d{\rm Vol_{M} }\nonumber
\end{align}
\begin{align}
&+\int_{\partial M}\Big[
\frac{1}{6} \pi  \left(10 \pi ^2-3\right)
\Big(\sum_{l=1}^{6}\langle\nabla_{J(e_{l})}^{L}(Je_{6}), e_{l}\rangle-\sum_{l=1}^{6}g^{M}\left(J(\frac{\partial}{\partial{x_{6}}}), \frac{\partial}{\partial{x_{6}}}\right)\langle\nabla_{J(e_{l})}^{L}(\frac{\partial}{\partial{x_{6}}}), e_{l}\rangle\Big)\nonumber\\
&-\frac{\pi }{8}\sum_{\beta=1}^{6}\langle J(e_{\beta}), e_{6}\rangle^{2}\Big(\sum_{l=1}^{6}\langle\nabla_{J(e_{l})}^{L}(Je_{6}), e_{l}\rangle-\sum_{l=1}^{6}g^{M}\left(J(\frac{\partial}{\partial{x_{6}}}), \frac{\partial}{\partial{x_{6}}}\right)\langle\nabla_{J(e_{l})}^{L}(\frac{\partial}{\partial{x_{6}}}), e_{l}\rangle\Big)\nonumber\\
&-\frac{\pi ^3}{3}\sum_{\beta=1}^{6}\sum_{i=1}^{5}\langle J(e_{\beta}), e_{i}\rangle^{2}\Big(\sum_{l=1}^{6}\langle\nabla_{J(e_{l})}^{L}(Je_{6}), e_{l}\rangle-\sum_{l=1}^{6}g^{M}\left(J(\frac{\partial}{\partial{x_{6}}}), \frac{\partial}{\partial{x_{6}}}\right) \langle\nabla_{J(e_{l})}^{L}(\frac{\partial}{\partial{x_{6}}}), e_{l}\rangle\Big)\nonumber\\
&-\pi ^3\Big(\sum_{l=1}^{6}g^{M}(J(e_{l}), (\nabla^{L}_{e_{l}}J)e_{6})
-\sum_{l=1}^{6}g^{M}(J(\frac{\partial}{\partial{x_{6}}}), (\nabla^{L}_{e_{l}}J)e_{l})+\sum_{l=1}^{6}g^{M}(J(e_{l}), (\nabla^{L}_{e_{6}}J)e_{l})\Big)\nonumber\\
&+\frac{\pi }{32}h'(0)\sum_{l,\beta=1}^{6}\langle J(e_{\beta}), e_{6}\rangle^{2} \langle J(e_{l}), e_{6}\rangle^{2}
+\frac{16 \pi ^3}{15}h'(0)\sum_{l,\beta=1}^{6}\sum_{i=1}^{5}\langle J(e_{l}), e_{6}\rangle^{2} \langle J(e_{\beta}), e_{i}\rangle^{2}\nonumber\\
&-\frac{3 \pi}{32}h'(0)\sum_{l=1}^{6}\sum_{i=1}^{5}\langle J(e_{i}), e_{6}\rangle^{2} \langle J(e_{l}), e_{6}\rangle^{2}
-\frac{5 \pi ^3}{12}h'(0)\sum_{l=1}^{6}\sum_{\nu,i=1}^{5}\langle J(e_{\nu}), e_{6}\rangle^{2} \langle J(e_{l}), e_{i}\rangle^{2}\nonumber\\
&-\frac{\pi ^3}{3}h'(0)\sum_{l=1}^{6}\sum_{\nu,i=1}^{5}\langle J(e_{l}), e_{6}\rangle^{2} \langle J(e_{\nu}), e_{i}\rangle^{2}
-\frac{3\pi ^3}{8}h'(0)\sum_{i=1}^{5}\langle J(e_{6}), e_{6}\rangle \langle J(e_{i}), e_{i}\rangle\nonumber\\
&-\frac{\pi ^3}{6}h'(0)\sum_{\nu,i=1}^{5}\langle J(e_{i}), e_{i}\rangle\langle J(e_{\nu}), e_{\nu}\rangle
+\frac{3}{80} \pi  \left(128 \pi ^2+5\right)h'(0)\sum_{l=1}^{6}\langle J(e_{l}), e_{6}\rangle^{2}\nonumber\\
&-\frac{5\pi ^3}{8}h'(0)\sum_{i=1}^{5}\langle J(e_{i}), e_{6}\rangle^{2}
+\frac{7 \pi ^3}{30}h'(0)\sum_{l=1}^{6}\sum_{i=1}^{5}\langle J(e_{l}), e_{i}\rangle^{2}\nonumber\\
&-\frac{7 \pi ^3}{10}h'(0)\sum_{l,\beta=1}^{6}\sum_{i=1}^{5}\langle J(e_{l}), e_{i}\rangle\ \langle J(e_{\beta}), e_{i}\rangle \langle J(e_{l}), e_{6}\rangle \langle J(e_{\beta}), e_{6}\rangle \nonumber\\
&+\frac{\pi ^3}{2}h'(0)\sum_{l=1}^{6}\sum_{\nu,i=1}^{5} \langle J(e_{\nu}), e_{i}\rangle \langle J(e_{l}), e_{i}\rangle \langle J(e_{\nu}), e_{6}\rangle \langle J(e_{l}), e_{6}\rangle \nonumber\\
&+\frac{\pi }{32}h'(0)\sum_{l=1}^{6}\sum_{i=1}^{5} \langle J(e_{l}), e_{6}\rangle^{2} \langle J(e_{i}), e_{i}\rangle \langle J(e_{6}), e_{6}\rangle\nonumber\\
&+\frac{\pi ^3}{4}h'(0)\sum_{l=1}^{6}\sum_{\nu,i=1}^{5}\langle J(e_{l}), e_{i}\rangle^{2} \langle J(e_{\nu}), e_{\nu}\rangle \langle J(e_{6}), e_{6}\rangle
\Big]\Omega_4d{\rm Vol_{\partial M}}.\nonumber
\end{align}
where $s$ is the scalar curvature.
\end{thm}

\section{Acknowledgements}

The author was supported in part by  NSFC No.11771070. The author thanks the referee for his (or her) careful reading and helpful comments.

\vskip 1 true cm


\bigskip
\bigskip

\noindent {\footnotesize {\it S. Liu} \\
{School of Mathematics and Statistics, Northeast Normal University, Changchun 130024, China}\\
{Email: liusy719@nenu.edu.cn}

\noindent {\footnotesize {\it Y. Wang} \\
{School of Mathematics and Statistics, Northeast Normal University, Changchun 130024, China}\\
{Email: wangy581@nenu.edu.cn}

\clearpage
\section*{Statements and Declarations}

Funding: This research was funded by National Natural Science Foundation of China: No.11771070.\\

Competing Interests: The authors have no relevant financial or non-financial interests to disclose.\\

Author Contributions: All authors contributed to the study conception and design. Material preparation, data collection and analysis were performed by Siyao Liu and Yong Wang. The first draft of the manuscript was written by Siyao Liu and all authors commented on previous versions of the manuscript. All authors read and approved the final manuscript.\\

Availability of Data and Material: The datasets supporting the conclusions of this article are included within the article and its additional files.\\

\end{document}